\documentclass[a4paper,12pt,oneside]{article}

\usepackage[utf8]{inputenc}
\usepackage[english]{babel} % Cambiado a ingl\'{e}s
\usepackage{hyperref}
\usepackage{authblk}
\usepackage{amsmath}   % mathematical environments
\usepackage{amsfonts}  % mathematical fonts
\usepackage{amssymb}   % additional symbols
\usepackage{amsthm}    % theorems and definitions
\usepackage{graphicx}
\usepackage{xcolor}
\usepackage{float}
\usepackage{subcaption}
\usepackage{geometry}
\theoremstyle{plain} % style with italics (for "strong" results)
\newtheorem{theorem}{Theorem}[section]
\newtheorem{lemma}[theorem]{Lemma}
\newtheorem{proposition}[theorem]{Proposition}
\newtheorem{corollary}[theorem]{Corollary}
\newtheorem{problem}[theorem]{Problem}
\theoremstyle{definition} % roman style (for definitions, examples)
\newtheorem{definition}[theorem]{Definition}

\theoremstyle{remark} % roman style (for comments or notes)
\newtheorem{remark}[theorem]{Remark}

\title{Backward Homogeneous Heat Equation with $\Psi$-Caputo Time Derivative and Spectral Fractional Laplacian}

\author[,1,2]{Ojeda, Nelson R.\thanks{\href{mailto:nojeda@santafe-conicet.gov.ar}{nojeda@santafe-conicet.gov.ar}}}
\author[,1,3]{Spies, Ruben D. \thanks{\href{mailto:rspies@santafe-conicet.gov.ar}{rspies@santafe-conicet.gov.ar}}}

\affil[1]{Instituto de Matemática Aplicada del Litoral, CONICET-UNL, Santa Fe, Argentina}

\affil[2]{Facultad de Recursos Naturales -Universidad Nacional de Formosa, Argentina}

\affil[3]{Departamento de Matemática, Facultad de Ingeniería Química, Universidad Nacional del Litoral, Santa Fe, Argentina}

\begin{document}

\maketitle
\begin{abstract}
    \textcolor{black}{In this work, we study a fractional diffusion model in time and space, formulated in an open and bounded domain,
    involving a Caputo fractional derivative of order $\mu\in(0,1)$ with respect to a time scale function $\Psi$, together with a spectral
    fractional Laplacian operator of order $s\in(0,1)$. The associated direct problem is analyzed, as well as the inverse problem of temporal
     reconstruction from final data, which is inherently ill-posed. Through a spectral representation, the instability of the problem and the
     effect of the choice of the function $\Psi$ on the amplification of spectral modes are examined.}
\end{abstract}



%35R30 — Inverse problems (el tema central del trabajo). 26A33 — Fractional derivatives and integrals. 35K05 — Heat equation. 35P10 — Eigenfunction expansions; spectral theory.47A75 — Eigenvalue problems.
\section{Introduction}

Anomalous diffusion processes constitute a broad family of stochastic phenomena
that deviate from the classical behavior of the standard diffusion equation. In
particular, \emph{ultra-slow diffusion} processes are characterized by a
logarithmic growth of the mean squared displacement, of the form
\begin{equation}\nonumber
\langle x^2(t) \rangle \sim (\log t)^{\beta}, \qquad \beta > 0,
\end{equation}
where $\langle x^2(t) \rangle$ denotes the mean square diffusion  as defined in \cite{KosztolowiczDutkiewicz2021a}.

Notable examples include \emph{Sinai diffusion} (for $\beta=4$), associated with the motion of particles in static random force fields, and certain models in
polymer physics where $\beta=4/3$ (\cite{Sinai1983}, \cite{SchiesselSokolovBlumen1997}, \cite{LiangEtAl2019}).

In this context, De Gregorio and Garra (\cite{DeGregorio2021}) introduced and analyzed fractional heat equations of Hadamard type, in which the temporal
fractional derivative is defined with respect to the function $\Psi(t)\doteq\ln (a+bt)$ for $a,\,b>0$. These equations naturally lead to ultra-slow type
diffusion laws, where the effective time evolution is described on a logarithmic scale.

Following this model, it is natural to consider a more general extension in which physical time $t$ is replaced by an increasing real function \(\Psi \in
C^1(0,\infty)\) such that its derivative $\Psi'(t) \neq 0$ in $(0,\infty)$. In such a case, following the works in \cite{KosztolowiczDutkiewicz2021a}
and\cite{KosztolowiczDutkiewicz2021b}, the $\alpha$-fractional $\Psi$-diffusion equation takes the form:
\begin{equation}\nonumber
{}_{}^{C}D_{a+;\Psi}^{\alpha} u(x,t) = D_{\alpha}\,\partial_x^2 u(x,t),
\end{equation}
where $D_{\alpha}$ is a subdiffusion coefficient measured in units of $m^2/t^{\alpha}$ and ${}_{}^{C}D_{a+;\Psi}^{\alpha}$ denotes the Caputo fractional
derivative in the temporal variable with respect to $\Psi$ of order $\alpha>0$, defined in \cite{Almeida2017}. By performing the change of variable $\tau =
\Psi(t)$, it is found that the above operator is equivalent to an ordinary fractional derivative in the variable $\tau$. Therefore, if in the $\tau$ scale  the
mean displacement satisfies $\langle x^2(\tau)\rangle \sim \tau^{\alpha}$, then, in the physical time $t$, there holds
\begin{equation}\nonumber
\langle x^2(t)\rangle \sim (\Psi(t))^{\alpha}.
\end{equation}
Thus this relationship shows that the choice of $\Psi$ determines the effective temporal law of the diffusion process.

In the case of fractional diffusion in the spatial variable, as an important background in this type of problems involving fractional and/or non-local
operators in space and time, a recent work by Li (\cite{Li2022}) studies the following non-homogeneous fractional equation:
\begin{equation}\label{eq.FBHEDeltaPsiR}
    \left\{\begin{aligned}
        {}_{}^{C}D_{a+;\Psi}^{\alpha} u(\mathbf{x},t) + (-\Delta)^s u(\mathbf{x},t) &= f(\mathbf{x},t),
        & \mathbf{x}\in \mathbb{R}^d,\,t>a  \\
        u(\mathbf{x},a^+) &= u_a(\mathbf{x}), & \mathbf{x}\in \mathbb{R}^d,
    \end{aligned}\right.
\end{equation}
where ${}_{}^{C}D_{a+;\Psi}^{\alpha}$ denotes the $\Psi$-Caputo derivative with $\alpha \in (0,1)$, the source term $f(\mathbf{x},t)$ and the initial data
$u_{a}(\mathbf{x})$ are given functions, and $(-\Delta)^{s}$ represents the $s$-fractional Laplacian with $s \in (0,1)$ given by the following integral
representation (see \cite{DiNezza2012}):

\begin{equation}\label{eq:LaplacianoINTR}
    (-\Delta)^{s} v(\mathbf{x})
    = C(d,s)\, \text{P.V.} \int_{\mathbb{R}^{d}}
    \frac{v(\mathbf{x})-v(\mathbf{y})}{|\mathbf{x}-\mathbf{y}|^{d+2s}}\, d\mathbf{y},
    \qquad \mathbf{x}\in \mathbb{R}^{d},
\end{equation}
where P.V. is the abbreviation for ``principal value" and $C(d,s)$ is a normalization constant. If $\mu=1$ and $s=1$ in \eqref{eq.FBHEDeltaPsiR}, it reduces to
the ordinary parabolic heat equation in $\mathbb{R}^d$. The general solution of \eqref{eq.FBHEDeltaPsiR} is obtained using Laplace transform method in the
temporal variable and Fourier transform in the spatial variable.

In the field of inverse problems, particularly in forward/backward heat equation type problems, there are already some interesting studies for specific
functions of $\Psi$. Indeed, taking $\Psi(t)=t$, leads to the classical Caputo fractional model associated with standard subdiffusion, first studied by Liu et
al. in 2010 (\cite{Liu2010}) in 2010, more recently by Nair et al. in 2024 (\cite{NairDanumjaya2024}) for the homogeneous case in one variable in the domain
$[0,\pi]$, and by Alavani et al. for the non-homogeneous case (\cite{AlavaniDanumjayaNair2024}). If we take $\Psi(t)=\log t$ we recovers the Hadamard-type
fractional model, corresponding to ultra-slow diffusion, studied in \cite{YangCaoLi2023}, \cite{LiLi2021a} and \cite{LiLi2021b,LiLi2021c}. Another interesting
case corresponds to the choice of $\Psi(t)\doteq t^\rho/\rho$, for $\rho>0$, which leads to the fractional model proposed by Katugampola, studied in
\cite{SinghAlshehriMomani2022}, \cite{BouchamaAriouaMerzougui2022} and \cite{OmabaSulaimani2022}.

%%%%%%%%%%%%%%%%%%%%%%%%%%%%%%%%%%%%%%%%%%%%%%%%%%%%%%%%%%%%%%%%%%%%%%%%%%%%%%%%%%%%%%%%%%%%%%%%%%%%%%%%%%%%%%%%%%% H  A  S  T  A

Note that problem \eqref{eq.FBHEDeltaPsiR} is studied in an unbounded domain
such as the entire $\mathbb{R}^d$. In this context, given an open and bounded
set $\Omega \subset \mathbb{R}^d$, we consider a homogeneous direct (forward)
problem formulated in a bounded cylindrical domain \(\Omega_{a,T}=\Omega \times
[a,T]\) for $0\leq a\leq t < T$, which incorporates a fractional derivative in
time of order \(\mu \in (0,1)\), in the Caputo sense with respect to the
function \(\Psi\), and a spatial diffusion of fractional order \(s \in (0,1)\),
modeled by the spectral fractional Laplacian operator. The problem is expressed
as:

\begin{equation}\label{eq.FBHEDeltaHOMOGPsi}
(PH_\Psi)\left\{\begin{aligned} {}_{}^{C}\partial^{\mu}_{a+;\Psi}
u(\mathbf{x},t)  &= -(-\Delta_{D,0})^s u(\mathbf{x},t),
& (\mathbf{x},t) \in \Omega_T, \\
u(\mathbf{x},t) &= 0, & (\mathbf{x},t) \in \partial\Omega \times (a,T], \\
u(\mathbf{x},a) &= u_a(\mathbf{x}), & \mathbf{x} \in \Omega,
\end{aligned}\right.
\end{equation}
where $\Delta_{D,0}$ denotes the Homogeneous Dirichlet spectral Laplacian, as defined in Definition \ref{Dirichlet_Laplacian}.

Similarly, the associated inverse problem is introduced:

\begin{quote}
$(IP)_{\Psi,\mu,s}:\quad$ \textit{Given the final state \(u_T \doteq u_{\mu,s}(\cdot,T)\) for some \(T>0\), determine the state \(u_{\mu,s}(\cdot,t)\) for all
\(0\leq a \le t < T\).}
\end{quote}

We are interested in studying the amplification of high-frequency modes with respect to the spectral growth of the eigenvalues $\lambda_k$ of the Laplacian. In
this sense, we will analyze whether the fractional backward $\Psi$-problem constitutes an example of an inverse problem whose degree of ill-posedness is lower
than that of classical problems. Furthermore, we will establish a stability hierarchy based on temporal evolution: we will demonstrate that, while the recovery
of the initial state ($t \to a^+$) remains a \textit{moderately ill-posed} problem of an algebraic nature, finding the solution at any intermediate time $t \in
(a, T)$ is a \textbf{well-posed} problem. This transition not only smooths the intrinsic instability of the diffusion process but also allows interpreting the
temporal flow itself as a natural regularization family for the inversion of operators involving the Mittag--Leffler function.

\section{Preliminaries}

\subsection{\texorpdfstring{$\Psi$}{Psi}-Fractional Calculus}

\subsubsection{Fractional Operators}

Several integral operators generalize the concept of an integral of fractional order $\mu$ (where $\mu>0$) of a function $f$:
\begin{eqnarray}
    \label{eq:IntRL}
    \text{ Riemann-Liouville }
    & \left({}^{RL}I_{a+}^{\mu} f\right)(t)
   \doteq\frac{1}{\Gamma(\mu)}
    \int_a^t (t-s)^{\mu-1} f(s)\, ds\\
    \label{eq:IntRLHada}
    \text{ Hadamard }
    & \left({}_{}^{H}I_{a+}^{\mu} f\right)(t)
   \doteq\frac{1}{\Gamma(\mu)}
    \int_a^t
    \left(\ln \frac{t}{s}\right)^{\mu-1}
    \frac{f(s)}{s}\, ds\\
    \label{eq:IntRLKatu}
    \text{ Katugampola }
    &  \left({}_{}^{K} I_{a+}^{\mu} f\right)(t)
   \doteq \frac{\rho^{1-\mu}}{\Gamma(\mu)}
    \int_a^t
    \bigl(t^\rho - s^\rho\bigr)^{\mu-1}
    f(s)s^{\rho-1}\, ds,\, \rho>0,
\end{eqnarray}
where $\Gamma(\cdot)$ denotes the classical gamma function. Fractional
derivatives are defined from fractional integrals and usually adopt a differintegral or Riemann-Liouville form $D^{\mu} \doteq D^{n} I^{\,n-\mu}$ and  the
integrodifferential or Caputo form $D^{\mu} \doteq I^{\,n-\mu}D^{n}$, where $I$ can
be any of the previously introduced fractional integrals and $D^{n}$ denotes
some integer-order differential operator. For example, if for $\mu >0$ we
denote $n = \lfloor \mu \rfloor + 1 \in \mathbb{N}$ if $\mu\not\in\mathbb{N}$
and $n=\mu$ if $\mu\in\mathbb{N}$, we have the following fractional operators:
\begin{eqnarray}
        \label{eq:DerRL}\nonumber
    &\text{ Riemann-Liouville }
    &\left({}^{RL}D_{a+}^{\mu} f\right)(t)
\doteq\left( \frac{d}{dt} \right)^{n}
{}^{RL}I_{a+}^{\,n-\mu} f(t)\\
    \label{eq:DerCAP}
&\text{ Caputo} &\left({}^{C}D_{a+}^{\mu} f\right)(t) \doteq
{}^{RL}I_{a+}^{\,n-\mu} \left( \frac{d}{dt} \right)^{n}f(t).
\end{eqnarray}

It should be noted that the operators defined in \eqref{eq:IntRL},
\eqref{eq:IntRLHada}, and \eqref{eq:IntRLKatu} can be unified through a
fractional integral defined with respect to a strictly increasing function
$\Psi$, as follows. 
\begin{definition}[{\cite{Podlubny1999,Fahad2021,Almeida2017, KilbasSrivastavaTrujillo2006, SamkoKilbasMarichev1993}}]
    Let $\mu$ be a real number such that $\mu > 0$, $-\infty < a < b < \infty$, $f$ an integrable function defined on $[a,b]$, and $\Psi \in C^1([a,b])$ an increasing function such that $\Psi'(t) \neq 0$ for all $t \in [a,b]$.
    Then, \textbf{the left-sided Riemann--Liouville $\mu$-fractional integral of $f$ with respect to $\Psi$}, is defined as
    \begin{equation}
        \label{eq:IntRLpsileft}
        \left(I_{a+;\Psi}^{\mu} f\right)(t)
        = \frac{1}{\Gamma(\mu)} \int_a^t \bigl(\Psi(t)-\Psi(s)\bigr)^{\mu-1}\Psi'(s) f(s)\, ds.
    \end{equation}
\end{definition}
Note that the choices $\Psi(t)=t$, $\Psi(t)=\ln t$ and
$\Psi(t)=\dfrac{t^\rho}{\rho}$ recover the Riemann--Liouville,
Hadamard, and Katugampola operators, respectively.

\begin{remark}
The term \textit{``left-sided operators"} is attributed to the fact that the lower
limit of the integral is considered a fixed point $a$. Thus, we can
analogously define \textbf{the right-sided Riemann--Liouville $\mu$-fractional
integral of $f$ with respect to $\Psi$} as
\begin{equation}\label{eq:IntRLpsi}
    \left(I_{b-;\Psi}^{\mu} f\right)(t)
    = \frac{1}{\Gamma(\mu)} \int_t^b \bigl(\Psi(s)-\Psi(t)\bigr)^{\mu-1}\Psi'(s) f(s)\, ds.
\end{equation}
\end{remark}

Likewise, we provide the following definition:
\begin{definition}[{\cite{Almeida2017}}]
Let $\mu\in \mathbb{R}$, $\mu > 0$, $-\infty < a < b <
\infty$, $n \doteq \lceil \mu \rceil$, and let $f, \Psi \in C^n([a,b])$ be functions such that $\Psi'(t) > 0$ for all $t \in [a,b]$. Then, the $\Psi$- Caputo  $\mu$-fractional derivative of $f$ is defined as
\begin{equation}\label{eq:DerCapsi}
\left({}_{}^{C}D_{a+;\Psi}^{\mu} f\right)(t) \doteq \left(I_{a+;\Psi}^{n-\mu}
f_{\Psi}^{[n]}\right)(t),
\end{equation}
where $f_{\Psi}^{[n]}(t)\doteq\left( \frac{1}{\Psi'(t)} \frac{d}{dt} \right)^n f(t)$
for all $t\in [a,b]$.
\end{definition}

In this work we are interested in studying the case  $\mu \in (0,1)$, in which case
 equation \eqref{eq:DerCapsi} takes de form
\begin{equation}\label{eq:DerCapsi1}
   \left({}_{}^{C}D_{a+;\Psi}^{\mu} f\right)(t)
= \frac{1}{\Gamma(\mu)} \int_a^t \bigl(\Psi(t)-\Psi(s)\bigr)^{-\mu} f'(s)\, ds.
\end{equation}
Note that if $\mu=1$ in equation \eqref{eq:DerCapsi1}, for a sufficiently smooth
$f$, we recover the first-order derivative, and 
$\left({}_{}^{C}D_{a+;\Psi}^{1} f\right)=\dfrac{1}{\Psi'(t)} \dfrac{d}{dt}$
holds.

Note that in the integrodifferential representation of the operator
\eqref{eq:DerCapsi1}, all the information from the classical derivative $f'(s)$
for $s \in [a,t]$ with $t\in [a,b]$ is used, which we call the \textit{``memory
effect"} of the fractional derivatives of the function $f$. It is known from the
fractional calculus literature that the essential difference between the
first-order derivative operator in the temporal variable $\partial_t$ and the
fractional operator $\partial_t^{\mu}$ for $\mu \in (0,1)$ occurs mainly when
performing a temporal inversion of the type $\tau = T - t$, where it is well
known that $\partial_t = -\partial_\tau$, whereas $\partial_t^{\mu} \neq
\partial_\tau^{\mu}$.

In the work of Fahad et al. (\cite{Fahad2021}), the authors revisited Samko's main idea that
the operators \eqref{eq:IntRLpsileft} and \eqref{eq:DerCapsi} can be
represented as the conjugation of the fractional operators \eqref{eq:IntRL} and
\eqref{eq:DerCAP} via the composition operation with $\Psi$ and $\Psi^{-1}$. More precisely, 

\begin{equation}\label{eq:psi_integral}
    {}_{}^{RL}I^{\mu}_{a+;\Psi}
    =Q_{\Psi}\circ {}_{}^{RL}I^{\mu}_{\Psi(a)+}\circ Q_{\Psi}^{-1},
\end{equation}
\begin{equation}
    \label{eq:psi_caputo}
    {}_{}^{C}D^{\mu}_{a+;\Psi}
    =
    Q_{\Psi}\circ {}_{}^{C}D^{\mu}_{\Psi(a)+}\circ Q_{\Psi}^{-1},
\end{equation}
where the linear operator $Q_{\Psi}$ is defined as

\begin{equation}\label{eq:def_Qpsi}
    (Q_{\Psi}f)(t)\doteq f(\Psi(t)), \quad t\in [a,b].
\end{equation}
This type of calculus is known as Operational Calculus and serves to deduce or
infer more general results involving the function $\Psi$ from classical results
of fractional calculus. 
Thus for instance, whe have the following existence
and uniqueness theorem:

\begin{theorem}\label{teo:PVPsiCap}
    For $0 < \mu \leq 1$, $0\leq a < t\leq b$, and $\lambda \in \mathbb{R}$, the fractional initial value problem 
    \begin{equation}\label{eq:PVIfracPsi}
    \left\{\begin{aligned}
        {}^{C}_{}D^{\mu}_{a;\Psi} y(t)  &=  \lambda y(t)\\
        y(a^+) &= y_a,
    \end{aligned}
    \right.
    \end{equation}
    has a unique solution given by
    \begin{equation}\label{eq:SolPVIfracPsi}
        y(t) = y_a\,E_{\mu}\!\big(\lambda(\Psi(t)-\Psi(a))^{\mu}\big),
    \end{equation}
    where $E_{\alpha}$ denotes the \textit{Mittag-Leffler function} defined by the hypergeometric series
\begin{equation}
E_{\alpha}(z) = \sum_{n=0}^{\infty} \frac{z^n}{\Gamma(n\alpha+1)}, \quad
z,\,\alpha \in \mathbb{C}.
\end{equation}
\end{theorem}
\begin{proof}
    Suppose that $\Psi \in C^{1}([a,b])$ and $\Psi'(t)>0$ on $(a,b)$. Then $\Psi$ is a bijection from $[a,b]$ onto $[\Psi(a),\Psi(b)]$ and $\Psi^{-1}\in C^{1}$ exists. In particular, the operator $Q_{\Psi}$  in \eqref{eq:def_Qpsi} is a linear isomorphism between the space of absolutely continuous functions $AC([\Psi(a),\Psi(b)])$ and the space
    %\begin{equation}\nonumber
    $
    X_\Psi\doteq\{y\in C([a,b]) : y=z\circ\Psi,\ z\in AC([\Psi(a),\Psi(b)])\},
    $
    %\end{equation}
    whose inverse is given by $(Q_\Psi^{-1}g)(s)=g(\Psi^{-1}(s))$ for all $s \in [\Psi(a),\Psi(b)]$. Now let $y$ be the solution to problem \eqref{eq:PVIfracPsi}. We define the auxiliary function
    \begin{equation}
        \label{eq:variable_auxZ}
        z \doteq Q_\Psi^{-1}y,
    \end{equation}
    that is $z(s)=y(\Psi^{-1}(s))$ for all $s \in [\Psi(a),\Psi(b)]$. Then, if we apply the operator $Q_\Psi^{-1}$ to the function ${}^{C}D^{\mu}_{a;\Psi} y(t)$ and use the conjugation identity \eqref{eq:psi_caputo}, we obtain
    \begin{equation}\nonumber
    Q_\Psi^{-1}\Big({}^{C}D^{\mu}_{a;\Psi} y\Big)
    =
    Q_\Psi^{-1}\Big(Q_\Psi\big({}^{C}D^{\mu}_{\Psi(a)+} z\big)\Big).
    \end{equation}
    Since $Q_\Psi^{-1}\circ Q_\Psi$ is the identity operator, it follows from (\ref{eq:PVIfracPsi}) that
    \begin{equation}\nonumber
    Q_\Psi^{-1}(\lambda y) =
  {}^{C}D^{\mu}_{\Psi(a)+} z.
    \end{equation}
    By the linearity of $Q_\Psi^{-1}$, we know that $Q_\Psi^{-1}(\lambda y) = \lambda Q_\Psi^{-1}y = \lambda z$. Therefore, for all $s \in [\Psi(a),\Psi(b)]$ there holds ${}^{C}D^{\mu}_{\Psi(a)+} z (s) = \lambda z(s)$. 
    
    Regarding the initial condition, we observe that $z(\Psi(a)^+) = y(a^+)$, so if $y(a^+)=y_a$, then $z(\Psi(a)^+)=y_a$.
    Thus, problem \eqref{eq:PVIfracPsi} is equivalent to solving the following fractional problem in the classical Caputo sense of order $\mu$,  $0<\mu\leq 1$:
    \begin{equation}
         \label{eq:PVIfrac} \left\{\begin{aligned} {}^{C}_{}D^{\mu}_{\Psi(a)+} z(s) &= \lambda z(s),\qquad s\in [\Psi(a),\Psi(b)],
            \\ z(\Psi(a)^+) &= y_a. \end{aligned} \right.
    \end{equation}
It is well known (see \cite{Podlubny1999,KilbasSrivastavaTrujillo2006, SamkoKilbasMarichev1993}) that (\ref{eq:PVIfrac}) has a unique solution given by
    \begin{equation}
        \label{eq:SolPVIfrac} z(s) = y_a\,E_{\mu}\!\big(\lambda (s-\Psi(a))^{\mu}\big),\qquad s\in [\Psi(a),\Psi(b)].
    \end{equation}
Therefore, taking $s=\Psi(t)$ for $t\in [a,b]$ in
\eqref{eq:variable_auxZ}, we get
\begin{equation}
     y(t)= y_a\,E_{\mu}\!\big(\lambda (\Psi(t)-\Psi(a))^{\mu}\big), \qquad t \in [a,b]
\end{equation}
which is what we wanted to prove.
\end{proof}

We conclude this section by recalling some properties and auxiliary results concerning the Mittag-Leffler function, which will be used in the following sections. The two-parameter Mittag-Leffler function is defined as
\begin{equation}
    \label{eq:Mitaggalphabeta}
    E_{\alpha,\,\beta}(z) \doteq \sum_{k=0}^{\infty} \frac{z^k}{\Gamma(\alpha k+\beta)},
    \quad z \in \mathbb{C},
\end{equation}
where $\alpha, \, \beta \in \mathbb{C}$. If $Re(\alpha)>0$, the series in
\eqref{eq:Mitaggalphabeta} converges for all $ z \in \mathbb{C}$ and defines an entire function of order $\frac{1}{Re(\alpha)}$ and type 1 (see \cite{Gorenflo2014}). Also, if $\beta=1$
in \eqref{eq:Mitaggalphabeta}, one has $E_{\alpha,\,1}=E_{\alpha}$.
Moreover, when $\alpha=\beta=1$, the Mittag-Leffler function reduces to the classical exponential function
\begin{equation}
    \nonumber
    E_{1,\,1}(z) = \sum_{k=0}^{\infty} \frac{z^k}{\Gamma(k+1)}=\sum_{k=0}^{\infty} \frac{z^k}{k!}=e^{z},
    \quad z \in \mathbb{C}.
\end{equation}
If $\alpha=0$ and $\beta=1$ in \eqref{eq:Mitaggalphabeta}
\begin{equation} \nonumber
    E_{0,\,1}(z) = \sum_{k=0}^{\infty} \frac{z^k}{\Gamma(0\cdot k+1)}=\sum_{k=0}^{\infty}z^k= \frac{1}{1-z},
    \quad |z|<1.
\end{equation}
If $\alpha =\frac{1}{2}$ and $\beta=1$ in \eqref{eq:Mitaggalphabeta}, then
\begin{equation}\nonumber
    E_{\frac{1}{2},\,1}(z) = \sum_{k=0}^{\infty} \frac{z^k}{\Gamma(\frac{1}{2}k+1)}=e^{z^2}\textrm{erfc}(-z),
    \quad z \in \mathbb{C}
\end{equation}
where $\textrm{erfc}(z)$ is the complementary error function $\textrm{erfc}(z)
\doteq \frac{2}{\sqrt{\pi}} \int_{z}^{\infty}  e^{-u^2}du$.

Finally, we recall several asymptotic properties of the Mittag-Leffler function.

\begin{theorem}
    Let $0 < \alpha < 2$, $\beta$ be an arbitrary complex number and $\mu$ be an arbitrary real number such that
     $\frac{\pi \alpha}{2} < \mu < \min\{\pi, \pi \alpha\}$. Then, for any integer $p \ge 1$, the following asymptotic expansion holds for all $z \in \mathbb{C}$ such that $\mu \le |\arg(z)| \le \pi$:
    \begin{equation}
        \label{eq:MittagAsin}
        E_{\alpha,\beta}(z)
        =
        -\sum_{k=1}^{p} \frac{z^{-k}}{\Gamma(\beta - \alpha k)}
        + \mathcal{O}\!\left(|z|^{-1-p}\right), \quad \textrm{ as } |z| \to \infty.
    \end{equation}
 Moreover, there exists a positive constant $C$ such that for all if $\beta \in \mathbb{R}$ there holds
    \begin{equation}
        \label{ineq:MittagAsin1}
        |E_{\alpha,\beta}(z)| \leq \frac{C}{1 + |z|}.
    \end{equation}
\end{theorem}
\begin{proof}
    The proof of the estimates \eqref{eq:MittagAsin} and \eqref{ineq:MittagAsin1} can be found in Theorem 1.4 and Theorem 1.6 of Podlubny (\cite{Podlubny1999}).
\end{proof}
\begin{corollary}
    \label{coro:MittagAsin}
    Suppose that $0 < \alpha_0 < \alpha_1 < 1$.
    Then there exist positive constants $C_1, C_2$, depending only on $\alpha_0$ and $\alpha_1$,
    such that for all $\alpha \in [\alpha_0, \alpha_1]$ there holds
    \begin{equation}
    \label{ineq:MittagAsin2}
    \frac{C_1}{\Gamma(1-\alpha)} \,\frac{1}{1+z}
    \;\leq\;
    E_{\alpha,1}(-z)
    \;\leq\;
    \frac{C_2}{\Gamma(1-\alpha)} \,\frac{1}{1+z},
    \qquad z \geq 0.
    \end{equation}
\end{corollary}

\begin{proof}
    See \cite{Liu2010}.
\end{proof}
\subsection{Fractional Laplacian}
\subsubsection{The Homogeneous Dirichlet Spectral Fractional Laplacian Operator $(-\Delta_{D,0})^s$}

Recall that given an open and bounded domain $\Omega \subset \mathbb{R}^d$ and the operator $-\Delta$, we can determine a complete system of eigenvalues and eigenfunctions $\{(\lambda_k, \varphi_k)\}_{k\in \mathbb{N}}$, which satisfies the homogeneous Dirichlet problem:
\begin{equation}
\label{eq:eigenlocal}
\begin{cases}
-\Delta u = \lambda\,u, & \text{in } \Omega,\\
u = 0, & \text{on } \partial\Omega.
\end{cases}
\end{equation}

It is well known in the literature that the Laplace operator with homogeneous
Dirichlet boundary conditions, denoted as $-\Delta_{D,0}$, is such that its
spectral representation for any function $u \in C_0^\infty(\Omega)$ is
\begin{equation}\nonumber
(-\Delta_{D,0})u = \sum_{k=1}^\infty \lambda_k\, u_k\, \varphi_k, \qquad
\text{where } \quad u_k = \int_\Omega
u(\mathbf{x})\,\varphi_k(\mathbf{x})\,d\mathbf{x}.
\end{equation}

Furthermore, since the resolvent of $-\Delta_{D,0}$ is compact in $L^2(\Omega)$,
its spectrum is discrete and consists solely of positive eigenvalues forming a
non-decreasing sequence
\begin{equation}\nonumber
0 < \lambda_1 \leq \lambda_2 \leq \cdots \leq \lambda_k \leq \cdots, \qquad
\text{with } \lim_{k\to\infty}\lambda_k = +\infty.
\end{equation}

The functions $\{\varphi_k : k\in \mathbb{N}\} \subset H_0^1(\Omega) \cap
C^\infty(\Omega)$ constitute an orthonormal basis of $L^2(\Omega)$. Note that
the scalars $u_k$ are the coefficients of
$u$ in the orthonormal basis $\{\varphi_k\}_{k=1}^\infty$ of $L^2(\Omega)$
associated with the operator $-\Delta_{D,0}$.

The following two definitions are based on the works \cite{Antil2018} and \cite{AntilBartels2017}.
\begin{definition}
    For all $s \geq 0$, we define the following fractional Sobolev space:
\begin{equation}\nonumber
\mathbb{H}^s(\Omega) \doteq \left\{ u = \sum_{k=1}^\infty u_k \varphi_k \in
L^2(\Omega) \,:\, \|u\|^2_{\mathbb{H}^s(\Omega)} \doteq \sum_{k=1}^\infty
\lambda_k^s u_k^2 < \infty \right\},
\end{equation}
where we recall that $\lambda_k$ are the eigenvalues of $-\Delta_{D,0}$, with
normalized eigenfunctions $\varphi_k$, and
\begin{equation}\nonumber
u_k = \langle u, \varphi_k\rangle_{L^2(\Omega)} = \int_\Omega u \varphi_k.
\end{equation}
Particularly, for $s=0$ we have $\mathbb{H}^0(\Omega)=L^2(\Omega)$, and for
$s=1$ we recover the classical Sobolev space
$\mathbb{H}^1(\Omega)=H^1(\Omega)$.
\end{definition}

\begin{definition} \label{Dirichlet_Laplacian}
Let $u \in C_c^\infty(\Omega)$, then we define the homogeneous Dirichlet
spectral Laplacian of $u$ of order $s$ with $0 < s \leq 1$ through the
following spectral representation:
\begin{equation}\nonumber
(-\Delta_{D,0})^s u \doteq \sum_{k=1}^\infty \lambda_k^s u_k \varphi_k, \qquad
\text{with } \quad u_k= \langle u, \varphi_k\rangle_{L^2(\Omega)}.
\end{equation}
\end{definition}

Note that for $s=1$, we have that $(-\Delta_{D,0})^1 = -\Delta_{D,0}$ in
$L^2(\Omega)$. Recall that by standard density techniques of
$C_c^\infty(\Omega)$, we can extend this definition to $\mathbb{H}^s(\Omega)$. Thus, if we
consider functions $u$ and  $v = \sum_{k=1}^\infty v_k \varphi_k \in
\mathbb{H}^s(\Omega)$, we obtain the following bound:
 \begin{equation}
    \int_\Omega (-\Delta_{D,0})^s u(\mathbf{x}) v(\mathbf{x}) d\mathbf{x}
\leq \|u\|_{\mathbb{H}^s(\Omega)} \, \|v\|_{\mathbb{H}^s(\Omega)}.
\end{equation}
In this way, the operator $(-\Delta_{D,0})^s$ extends to an operator 
from $\mathbb{H}^s(\Omega)$ to its dual $\mathbb{H}^{-s}(\Omega)$ by density.
Furthermore, we note that in this case
\begin{equation}\nonumber
\|u\|_{\mathbb{H}^s(\Omega)} = \|(-\Delta_{D,0})^{s/2} u\|_{L^2(\Omega)}.
\end{equation}

\begin{remark}\label{obs:Servadei}
    In 2014, Servadei and Valdinoci \cite{ServadeiValdinoci2014} studied the fractional Laplacian defined in a different way equivalent to \eqref{eq:LaplacianoINTR}. They pointed out that if we compare the operators $(-\Delta)^{s} v(\mathbf{x})$ and $(-\Delta_{D,0})^{s} v(\mathbf{x})$, although in part of the literature both operators are often used interchangeably, they are actually different. Indeed, the eigenvalue problem associated with a \textit{non-local} homogeneous Dirichlet boundary condition for the operator $(-\Delta)^{s} v(\mathbf{x})$ is
    \begin{equation}\label{eq:s_prob_int}
        \begin{cases}
            (-\Delta)^s u = \lambda u & \text{in } \Omega,\\
            u = 0  & \text{in } \mathbb{R}^n \setminus \Omega.
        \end{cases}
    \end{equation}
\noindent
    They proved that the operators $(-\Delta)^s$ and $(-\Delta_{D,0})^s$ are not the same in $\Omega$, as they possess different eigenvalues and eigenfunctions (under the  Dirichlet boundary conditions \eqref{eq:eigenlocal} and \eqref{eq:s_prob_int} respectively), given that the first eigenvalue of $(-\Delta)^s$ is strictly smaller than that of $(-\Delta_{D,0})^s$, and the eigenfunctions of $(-\Delta)^s$ are generally only H\"{o}lder continuous up to the boundary, in contrast to those of $(-\Delta_{D,0})^s$, which are as smooth up to the boundary as the domain $\Omega$ allows it.
\end{remark}

\section{Main Results}

From Remark \ref{obs:Servadei} it follows that the  study of diffusion problems in bounded domains using the
Laplacian $(-\Delta)^s$ will no yield the same results as studying them with the spectrally defined Laplacian $(-\Delta_{D,0})^s$. In what follows we plan to tackle the study of problem \eqref{eq.FBHEDeltaPsiR} in a bounded domain $\Omega\subset \mathbb{R}^d$.

\subsection{Forward Problem in a Bounded Domain}
Let $\Omega \subset \mathbb{R}^d$ be an open and bounded set. For $T
> a \ge 0$ we define the following set:
\begin{equation}\nonumber
    \Omega_{a,T} \doteq \Omega \times (a,T].
\end{equation}

We consider the following problem for the bounded set $\Omega_{a,T} \subset
\mathbb{R}^{d+1}$ with a fractional derivative in the temporal variable of
order $\mu \in (0,1)$ in the Caputo sense with respect to a function $\Psi$
(denoted by ${}_{}^{C}\partial^{\mu}_{a+;\Psi(t)}$), and of order $s \in (0,1)$
in the space variable, in the sense of the spectral fractional Laplacian (Definition \ref{Dirichlet_Laplacian}), that is:
\begin{equation}
    \label{eq.FBHEDeltaPsi}
(PFH_\Psi)  \left\{\begin{aligned}
        {}_{}^{C}\partial^{\mu}_{a+;\Psi(t)} u(\mathbf{x},t)  &= -(-\Delta_{D,0})^s u(\mathbf{x},t),
        & (\mathbf{x},t) \in \Omega_{a,T}, \\
        u(\mathbf{x},t) &= 0, &(\mathbf{x},t) \in \partial\Omega \times (a,T], \\
        u(\mathbf{x},a^+) &= u_a(\mathbf{x}), & \mathbf{x} \in \Omega,
    \end{aligned}\right.
\end{equation}
where $u_a(\mathbf{x})$ is the initial state satisfying the compatibility
condition
\begin{equation}\nonumber
    u_a(\mathbf{x}) = 0, \quad \mathbf{x} \in \partial\Omega.
\end{equation}
For $\Psi : [a,T] \to \mathbb{R}$, we
choose a strictly increasing function $\Psi \in C^1([a,T])$ and associate an
auxiliary function
\begin{equation}
    \label{eq:Psi_aux}
    \Psi_a(t) \doteq \Psi(t) - \Psi(a), \qquad t \in [a,T].
\end{equation}

We propose that the solution of ($PFH_\Psi$) admits the following series representation:
\begin{equation}
    \label{eq:representation_u}
    u(\mathbf{x},t) = \sum_{k=1}^{\infty} u_k(t) \,\varphi_k(\mathbf{x}),\qquad (\mathbf{x},t) \in \Omega_{a,T}
\end{equation}
where for each $k \in \mathbb{N}$, a priori, the functions $u_k(t)$ may depend
on $\mu$ and $s$. 

We have the following result.

\begin{theorem}[\textbf{Existence, Uniqueness and the Solution Operator}]
    Let $s, \mu \in (0,1)$, $0 \leq a < t \leq T$, the initial data function $u_a \in L^2(\Omega)$, and  $\Psi_a$ as in \eqref{eq:Psi_aux}.
  Let also the operator ${_a}\mathcal{K}_{\Psi(t)}^{s,\mu} : L^2(\Omega) \to L^2(\Omega)$ be defined by
    \begin{equation}
        \label{eq:AalphaPsi}
        {_a}\mathcal{K}_{\Psi(t)}^{s,\mu}(v)
        \doteq
        \sum_{k=1}^{\infty}
        E_{\mu}\!\big(-\lambda_k^{\,s}(\Psi_a(t))^{\mu}\big)\,v_k\,\varphi_k,
        \qquad v \in L^2(\Omega),
    \end{equation}
    where $\{(\lambda_k, \varphi_k)\}_{k\in \mathbb{N}}$ is the system of eigenvalues and eigenfunctions that solves problem \eqref{eq:eigenlocal}. Then $(PFH_\Psi)$ has a unique solution that can be writen in the form 
    \begin{equation}
        \label{eq:SolFBHEDeltaPsi}
        u_{\mu,s}(\mathbf{x},t)
        =
        {_a}\mathcal{K}_{\Psi(t)}^{s,\mu}(u_a)(\mathbf{x}), \qquad \text{for all } (\mathbf{x},t) \in \Omega \times [a,T].
    \end{equation}
\end{theorem}
\begin{proof}
    Let $u$ be a function given by \eqref{eq:representation_u}. Substituting into \eqref{eq.FBHEDeltaPsi}, we obtain:
    \begin{equation} \nonumber
       %{}_{a}^{C}\partial^{\mu}_{\Psi(t)} u(\mathbf{x},t)  &= -(-\Delta_{D,0})^s u(\mathbf{x},t),
        %\qquad (\mathbf{x},t) \in \Omega_{a,T}, \\
        \sum_{k=1}^{\infty} \left( {}^{C}_{}D^{\mu}_{a+;\Psi} u_k\right)(t) \,\varphi_k(\mathbf{x}) = - \sum_{k=1}^{\infty} \left( \lambda_k^{s} u_k(t)\right) \varphi_k(\mathbf{x}),
        \end{equation}
        that is
    \begin{equation} \label{zzz}
        \sum_{k=1}^{\infty} \left( \left( {}^{C}_{}D^{\mu}_{a+;\Psi} u_k\right)(t) + \lambda_k^{s} u_k(t)\right) \varphi_k(\mathbf{x}) = 0.
    \end{equation}
\noindent
    Since the family of functions $\{\varphi_k : k \in \mathbb{N}\}$ is an orthonormal basis of $L^2(\Omega)$, it follows that all coefficients in (\ref{zzz}) must vanish, which leads to the following initial value problems in the Caputo sense, decoupled in the temporal variable with respect to $\Psi$, of order $0 < \mu < 1$:
    \begin{equation}\label{eq:fraccPsi}
        \left\{\begin{aligned}
             \left({}^{C}_{}D^{\mu}_{a+;\Psi} u_k\right)(t)  &= -\lambda_k^{s} u_k(t),
            && k \in \mathbb{N}, \; t \in (a,T], \\
            u_k(a) &= (u_a)_k,
        \end{aligned}\right.
    \end{equation}
where $(u_a)_k=\langle u_a, \varphi_k\rangle$. 
    Problem \eqref{eq:fraccPsi} is now solved using Theorem \ref{teo:PVPsiCap}, which guarantees the existence and uniqueness of the solution for the fractional problem with the Caputo derivative with respect to $\Psi$. For each $k\in\mathbb{N}$, the solution is given by:
    \begin{equation}
    \label{eq:coef_representation_u}
    u_k(t) = E_{\mu}\!\big(-\lambda_k^s(\Psi_a(t))^{\mu}\big)\,(u_a)_k,\quad t \in [a,T].
    \end{equation}
    Finally, substituting the coefficients obtained in \eqref{eq:coef_representation_u} into the series expansion \eqref{eq:representation_u}, we obtain the following expression:
    \begin{equation}
        \label{eq:solu}
        u_{\mu,s}(\mathbf{x},t)=\sum_{k=1}^{\infty} E_{\mu}\big(-\lambda_k^s\Psi_a(t))^{\mu}\big)\,(u_a)_k\, \varphi_k(\mathbf{x}), \qquad \text{for all } (\mathbf{x},t)\in\Omega\times [a,T].
    \end{equation}
    Note that if $u_a \in L^2(\Omega)$, then for each $t\in [a,T]$, $u_{\mu,s}(\cdot,t) \in L^2(\Omega)$. Indeed:
    \begin{align*}
        \|u_{\mu,s}(\cdot,t)\|^2_{L^2(\Omega)}&= \sum_{k=1}^{\infty} |E_{\mu}\big(-\lambda_k^s\Psi_a(t))^{\mu}\big)|^2\,(u_a)^2_k \\
        &\leq \sum_{k=1}^{\infty} ||E_{\mu}\big(-\lambda_k^s\Psi_a(\cdot))^{\mu}\big)||_{L^{\infty}(a,T)}^2\,(u_a)^2_k \qquad\qquad (cf. \eqref{eq:cotaMitagg1})\\
        &= \sum_{k=1}^{\infty} (u_a)^2_k \;=\;\|u_a\|^2_{L^2(\Omega)}<\infty,
    \end{align*}
where we have used that $||E_{\mu}\big(-\lambda_k^s\Psi_a(\cdot))^{\mu}\big)||_{L^{\infty}(a,T)}=1$, a fact that shall be proved below (see \eqref{eq:cotaMitagg1}). 
Hence, the solution operator ${_a}\mathcal{K}_{\Psi(t)}^{s,\mu}$ given by \eqref{eq:AalphaPsi} is bounded from $L^2(\Omega)$ to $L^2(\Omega)$.
\end{proof}

\begin{remark}[Particular cases of the operator ${_a}\mathcal{K}_{\Psi(t)}^{s,\mu}$]
    Some particular choices of the temporal transformation function $\Psi$ and the initial time parameter $a$,  allow us to recover well-known fractional operators. In fact:

    \begin{itemize}
        \item \textit{Caputo-Katugampola:} $\Psi(t)=\dfrac{t^\rho}{\rho}$, $0<\rho\leq 1$, and $a=0$. In this case
    \begin{equation}\label{eq.SolFBHEDeltaCaputokatumgampola}
            \nonumber
            {_0}\mathcal{K}_{\frac{t^\rho}{\rho}}^{s,\mu}(v) = \sum_{k=1}^\infty
            E_{\mu}\!\left(-\frac{\lambda_k^s\,t^{\rho\mu}}{\rho^\mu}\right)\,v_k.
        \end{equation}
        
        \item \textit{Classical Caputo:} $\Psi(t)=t$ and $a=0$. In this case
        \begin{equation}\label{eq.SolFBHEDeltaCaputo1}
            {_0}\mathcal{K}_{t}^{s,\mu}(v) = \sum_{k=1}^\infty
            E_{\mu}\!\big(-\lambda_k^st^{\mu}\big)\,v_k.
        \end{equation}
        If $\mu=s = 1$ in \eqref{eq.SolFBHEDeltaCaputo1}, we recover the classical heat semigroup $\{e^{(-\Delta_{D,0}) t}\}_{t\geq 0}$:
        \begin{equation}\label{eq.SolFBHEDeltaCaputo3}
            {_0}\mathcal{K}_{t}^{1,1}(v) = \sum_{k=1}^\infty
            e^{-\lambda_kt}\,v_k\,
            \varphi_k = e^{-\Delta t } v.
        \end{equation}
        \item \textit{Caputo-Hadamard:} $\Psi(t)=\ln t$ and $a>0$. Here
        \begin{equation}\label{eq.SolFBHEDeltaHadamar}
            \nonumber
            {_a}\mathcal{K}_{\ln t}^{s,\mu}(v)  = \sum_{k=1}^\infty
            E_{\mu}\!\left(-\lambda_k^s\left(\ln( t/a)\right)^{\mu}\right)\,v_k .
        \end{equation}
    \end{itemize}
\end{remark}

\subsection{On the Mittag-Leffler Function}

Note that the argument in the Mittag-Leffler function appearing in the solution operator
\eqref{eq:AalphaPsi}, for $0<s<1$ and $0<\mu<1$, has the form:
\begin{equation}
    \label{eq:argumentoZk}
    z_k = -\lambda_k^s \,[\Psi_a(t)]^\mu, \qquad k \in \mathbb{N},\; t\in [a,T].
\end{equation}

\begin{figure}[H]
    \centering
    \includegraphics[scale=0.5]{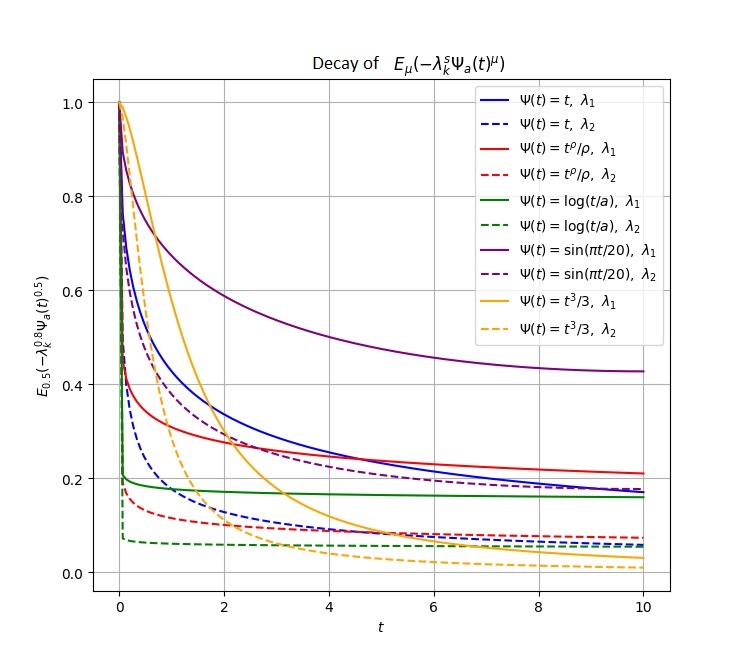}
    \caption{\scriptsize
        Comparison of the decay of the Mittag-Leffler function $E_{\mu}(-\lambda_k^{s}\Psi(t)^{\mu})$ for $\mu=0.5$ and $s=0.8$, considering two spectral values $\lambda_1=1$ (solid lines) and $\lambda_2=4$ (dashed lines). The curves correspond to different choices of the temporal scale function $\Psi(t)$ at final time $T=10$: $\Psi(t)=t$ with $a=0$, $\Psi(t)=t^{\rho}/\rho$ with $\rho=0.4$ and $a=0$, $\Psi(t)=\log(t/a)$ with $a=10^{-4}$, $\Psi(t)=\sin(\pi t/20)$ with $a=0$, and $\Psi(t)=t^3/3$ with $a=0$. The figure illustrates how the spectral parameter $\lambda_k$ and the choice of the function $\Psi(t)$ influence the decay rate of the relaxation process governed by the Mittag-Leffler function.
    }
        \label{fig:Mitagg}
\end{figure}

We are interested in studying the monotonicity and asymptotic behavior of the Mittag-Leffler function with respect to time and the eigenvalues, once the parameters $s$ and
$\mu$ are fixed. Figure \ref{fig:Mitagg} depicts the behavior of $E_\mu$ for some particular choices of the funcion $\Psi$.

The following result proves both temporal and spectral monotonicity of the Mittag-Leffler function.   
%%%%%%%%% DESDE -----------------------
\begin{lemma}
    \label{lem:Mittag}
    Let the parameters $s,\mu\in(0,1)$ and $0\leq a< T$ be fixed, and let $\Psi_a$ be as in \eqref{eq:Psi_aux}. Then, the following properties hold:
    \begin{itemize}
        \item[a)] \textbf{Spectral monotonicity}: For each fixed $t\in (a,T]$, the function defined by the mapping $\lambda \mapsto E_{\mu}\!\big(-\lambda \,[\Psi_a(t)]^\mu\big)$ is positive and strictly decreasing.
        \item[b)] \textbf{Temporal monotonicity}: For fixed $\lambda >0$, the mapping $t\mapsto E_{\mu}\!\big(-\lambda\,[\Psi_a(t)]^\mu\big)$ is positive, strictly decreasing, and satisfies:
        \begin{equation}
        \label{eq:cotaMitagg1}
        \Big\|
        E_{\mu}\!\big(-\lambda[\Psi_a(\cdot)]^\mu\big)
        \Big\|_{L^\infty(a,T)}
        =1.
        \end{equation}
        \item[c)] $\Psi$-\textbf{decay}: For any $\mu_0$, $\mu_1\in (0,1)$, $\mu_0<\mu_1$, there exist positive constants $C_3,C_4$ (depending on $\mu_0$ and $\mu_1$) such that for every $t\in(a,T]$ and $\mu\in[\mu_0,\mu_1]$ there holds:
        \begin{equation}
        \label{eq:cotaMitagg2}
        \frac{C_3}{1+\lambda \,[\Psi_a(t)]^\mu}
        \;\leq\;
        E_{\mu}\!\big(-\lambda \,[\Psi_a(t)]^\mu\big)
        \;\leq\;
        \frac{C_4}{1+\lambda \,[\Psi_a(t)]^\mu}.
        \end{equation}
        \item[d)] \textbf{Convergence}: For each fixed \(t\in(a,T]\), it holds that:
        \begin{equation}
        \label{eq:convergenceMitagg}
        E_{\mu}\!\big(-\lambda[\Psi_a(t)]^\mu\big)
        \overset{\lambda\to\infty}{\longrightarrow} \;0
        \end{equation}
    \end{itemize}
\end{lemma}

\begin{proof}
Properties a) and b) are deduced by considering \eqref{eq:argumentoZk} and the
fact that the mapping $z \mapsto E_\mu(-z)$ is completely monotonic for $z\geq
0$ given a fixed $\mu \in (0,1)$ (see \cite[Chapter~3]{Gorenflo2014}). Since
the Mittag-Leffler function reaches its maximum when the argument vanishes, where $E_{\mu}(0)=1$, \eqref{eq:cotaMitagg1} follows immediately.

For c), using Corollary \ref{coro:MittagAsin} with the argument given by \eqref{eq:argumentoZk} there follows: 
    \begin{equation}\nonumber
\frac{C_1}{\Gamma(1-\mu)} \,\frac{1}{1+\lambda[\Psi_a(t)]^\mu} \;\leq\;
E_{\mu}(-\lambda[\Psi_a(t)]^\mu) \;\leq\; \frac{C_2}{\Gamma(1-\mu)}
\,\frac{1}{1+\lambda[\Psi_a(t)]^\mu},
\end{equation}
where $C_1$ and $C_2$ are constants depending only on $\mu_0$ and $\mu_1$.  Letting $C_3\doteq \dfrac{C_1}{\Gamma(1-\mu_1)}$ and $C_4\doteq \dfrac{C_2}{\Gamma(1-\mu_0)}$ we obtain the bound \eqref{eq:cotaMitagg2}.

From the above results we then have that for every $t\in[a,T]$ and 
$\lambda >0$ there holds
    \begin{equation}\label{eq:cotakappa}
        0
        \le
            E_{\mu,1}\!\big(-\lambda[\Psi_a(t)]^\mu\big)
        \le
        \frac{C_4}{1+\lambda[\Psi_a(t)]^\mu}.
    \end{equation}
    Letting $\lambda \to +\infty$ in (\ref{eq:cotakappa}) yields the convergence \eqref{eq:convergenceMitagg}.
\end{proof}
\begin{remark}
It is well known that if $\{\lambda_k\}_{k\ge1}$ is the sequence of 
eigenvalues of the Laplacian in a bounded domain $\Omega \subset \mathbb{R}^d$
with homogeneous boundary conditions, Weyl's law (\cite{Grebenkov-2013}) implies that there exist
positive constants $c_1,c_2$ (depending only on $\Omega$) and 
$k_0\in\mathbb{N}$ such that:
\begin{equation}
    \label{eq:Weyl}
    c_1 k^{2/d}
    \;\le\;
    \lambda_k
    \;\le\;
    c_2 k^{2/d},
    \qquad \forall k \ge k_0.
\end{equation}
 Since $s>0$, the function $x\mapsto x^s$ is monotonically increasing on $(0,
\infty)$. Then \eqref{eq:Weyl} yields the following spectral $s$-Weyl estimate:
\begin{equation}
    \label{eq:Weyl_frac}
    c_1^{\,s} k^{2s/d}
    \;\le\;
    \lambda_k^{s}
    \;\le\;
    c_2^{\,s} k^{2s/d},
    \qquad \forall  k \ge k_0.
\end{equation}
\end{remark}
The following result shows the spectral asymptotic behavior of
the Mittag-Leffler function.

\begin{theorem}[Spectral asymptotic estimation of the Mittag-Leffler funtion]
    \label{coro:MittagAsin2}

    Let $\Omega\subset\mathbb{R}^d$ be a bounded domain and suppose $0<s<1$, $0<\mu_0<\mu_1<1$. Then there exist positive constants $m,M>0$ (depending only on $\mu_0,\mu_1,s$ and on the domain $\Omega$) such that for all $\mu\in[\mu_0,\mu_1]$ and $t\in(a,T]$ there exist $k_0\in\mathbb{N}$ (depending on $\mu$, $t$, $a$ and $s$) for which there holds:
    \begin{equation}\nonumber
    \frac{m}{k^{2s/d}[\Psi_a(t)]^\mu}
    \;\le\;
    E_\mu\!\big(-\lambda_k^{s}[\Psi_a(t)]^\mu\big)
    \;\le\;
    \frac{M}{k^{2s/d}[\Psi_a(t)]^\mu} \quad \forall k\ge k_0.
    \end{equation}
    In particular, the spectral modes of the solution operator (\ref{eq:AalphaPsi}), exhibit an asymptotic polynomial behavior of order $2s/d$, modulated by the temporal factor $[\Psi_a(t)]^{-\mu}$.
\end{theorem}
\begin{proof}
    From the asymptotic estimation of the Mittag-Leffler function (estimate \eqref{eq:cotaMitagg2} from Lemma \ref{lem:Mittag}), there exist positive constants $C_3,C_4$, depending only on $\mu_0$ and $\mu_1$, such that:
    \begin{equation}
        \label{eq:ML-lambda}
        \frac{C_3}{1+\lambda_k^{s}[\Psi_a(t)]^{\mu}}
        \;\le\;
        E_\mu\!\big(-\lambda_k^{s}[\Psi_a(t)]^{\mu}\big)
        \;\le\;
        \frac{C_4}{1+\lambda_k^{s}[\Psi_a(t)]^{\mu}}.
    \end{equation}
    Since $[\Psi_a(t)]^\mu>0$ and $\lambda_k\to+\infty$ as $k\to\infty$, there exists $k_0\in\mathbb{N}$ such that:
    \begin{equation}\nonumber
    \lambda_k^{s}[\Psi_a(t)]^\mu \ge 1,
    \qquad \forall k\ge k_0.
    \end{equation}
    Note that $k_0$ depends on $\mu$, $t$, $a$ and $s$.
    With these estimates into \eqref{eq:ML-lambda}, we obtain:
    \begin{equation}\nonumber
    \frac{C_3}{2\,\lambda_k^{s}[\Psi_a(t)]^{\mu}}
    \le
    E_\mu\!\big(-\lambda_k^{s}[\Psi_a(t)]^{\mu}\big)
    \le
    \frac{C_4}{\lambda_k^{s}[\Psi_a(t)]^{\mu}},
    \qquad \forall k\ge k_0.
    \end{equation}
    Finally, using the spectral $s$-Weyl estimate \eqref{eq:Weyl_frac}, we deduce that there exist positive constants $m\doteq\frac{C_3}{2c^s_2}$ and $M\doteq\frac{C_4}{c^s_1}$, independent of $k$, such that:
    \begin{equation}\nonumber
    \frac{m}{k^{2s/d}[\Psi_a(t)]^\mu}
    \le
    E_\mu\!\big(-\lambda_k^{s}[\Psi_a(t)]^\mu\big)
    \le
    \frac{M}{k^{2s/d}[\Psi_a(t)]^\mu},
    \qquad \forall k\ge k_0.
    \end{equation}
    This concludes the proof.
\end{proof}

\subsection{Properties of the Solution Operator}

In this section we present a few results related to some properties of the operator ${_a}\mathcal{K}^{s,\mu}_{\Psi(t)}$. In particular we characterize its spectrum and prove that it is self-adjoint, positive, compact, and injective.

\begin{proposition}
    \label{prop:Operator_Sol}
    Let ${_a}\mathcal{K}^{s,\mu}_{\Psi(t)}$ be the solution operator defined in\eqref{eq:AalphaPsi}. Then, for each $t\in(a,T]$, the following properties hold:
    \begin{enumerate}
        \item The point spectrum of the operator ${_a}\mathcal{K}^{s,\mu}_{\Psi(t)}$ is given by: 
        \begin{equation}\nonumber
        \sigma_p({_a}\mathcal{K}^{s,\mu}_{\Psi(t)})
        =
        \Big\{
        E_{\mu}\!\big(-\lambda_k^{\,s}(\Psi_a(t))^{\mu}\big)
        :\, k\in\mathbb{N}
        \Big\}
        \subset (0,\infty).
        \end{equation}
and  the family $\{\varphi_k\}_{k\in\mathbb{N}}$ constitutes the set of the corresponding eigenfunctions.
        \item The sequence of eigenvalues of ${_a}\mathcal{K}^{s,\mu}_{\Psi(t)}$ is positive and monotonically decreasing with respect to both the sequence of eigenvalues $\lambda^s_k$ and the variable $t\in(a,T]$.
        \item The operator ${_a}\mathcal{K}^{s,\mu}_{\Psi(t)}$ is self-adjoint, positive, compact, and injective.
    \end{enumerate}
\end{proposition}

\begin{proof}
For any $t\in (a,T]$ and $k\in\mathbb{N}$, we define
\begin{equation} \label{sigmak}
\sigma_k(t)\doteq E_{\mu}\!\big(-\lambda_k^{s}(\Psi_a(t))^{\mu}\big).
\end{equation}
 Note that \(\sigma_k(t)\in\mathbb{R}\) for all \(k\in
\mathbb{N}\). Let us recall that $\{\varphi_k\}_{k\in\mathbb{N}}$ is an orthonormal
basis of $L^2(\Omega)$, and by virtue of (\ref{eq:AalphaPsi}) the action of ${_a}\mathcal{K}^{s,\mu}_{\Psi(t)}$ on this basis is:
\begin{equation}\nonumber
{_a}\mathcal{K}^{s,\mu}_{\Psi(t)}(\varphi_k)=\sigma_k(t)\,\varphi_k,\qquad
k\in\mathbb{N},
\end{equation}
and consequently, its spectrum is pointwise being $\sigma_k(t)$ the eigenvalues.

Let \(u=\sum_{k} u_k\varphi_k\) and \(v=\sum_{k} v_k\varphi_k\) be two
functions  in \(L^2(\Omega)\). Then,
\begin{equation}\nonumber
\langle{_a}\mathcal{K}^{s,\mu}_{\Psi(t)} u, v\rangle = \sum_{k=1}^{\infty} (\sigma_k(t)\,
u_k)\,  \overline{v_k} = \sum_{k=1}^{\infty} u_k (\overline{\sigma_k(t)\,
v_k})  = \langle
u,{_a}\mathcal{K}^{s,\mu}_{\Psi(t)} v\rangle.
\end{equation}
Hence ${_a}\mathcal{K}^{s,\mu}_{\Psi(t)}$ is self-adjoint.

For positivity observe that for any $u \in
L^2(\Omega)$
\begin{equation}\nonumber
\langle {_a}\mathcal{K}^{s,\mu}_{\Psi(t)} u,u\rangle = \sum_{k=1}^{\infty} \sigma_k(t)
|u_k|^2 \ge 0,
\end{equation}
since \(\sigma_k(t)>0\) for all $t\in (a,T]$. Furthermore ${_a}\mathcal{K}^{s,\mu}_{\Psi(t)}$is strictly positive since  $\langle {_a}\mathcal{K}^{s,\mu}_{\Psi(t)} u,u\rangle =0$ if and only if $u=0$.

For the inyectivity, suppose that ${_a}\mathcal{K}^{s,\mu}_{\Psi(t)}\,v=0$ for a certain $v\in
L^2(\Omega)$. This clearly implies $\sigma_k v_k=0$ for all $k\in\mathbb{N}$. Now, for $t\in(a,T]$, we
have $\Psi_a(t)> 0$, and for $z\ge0$ and $0<\mu<1$, the
Mittag-Leffler function satisfies $E_\mu(-z)>0$. Therefore, $\sigma_k(t)>0$ for
all $k$. Hence $v_k=0$ for all $k$ and therefore
\begin{equation}\nonumber
\ker\!\big({_a}\mathcal{K}^{s,\mu}_{\Psi(t)}\big)=\{0\},
\end{equation}
which proves that the operator is injective.

Let us now show that ${_a}\mathcal{K}^{s,\mu}_{\Psi(t)}$ is compact. For that let $v=\sum_{k=1}^\infty
v_k\varphi_k \in L^2(\Omega)$ and consider the finite-rank approximations:
\begin{equation}\nonumber
K_N\,v\doteq\sum_{k=1}^N \sigma_k(t)\,v_k\,\varphi_k,\qquad N\in\mathbb{N}.
\end{equation}
Then \(K_N\) is of finite rank by construction. We will show that
\(\|{_a}\mathcal{K}^{s,\mu}_{\Psi(t)}-K_N\|_{\mathcal{L}(L^2)}\to 0\) as \(N\to\infty\). Now, the $\mathcal{L}(L^2)$-norm of the operator ${_a}\mathcal{K}^{s,\mu}_{\Psi(t)}-K_N$ coincides with its spectral radius, that is: 
\begin{equation}\nonumber
\|{_a}\mathcal{K}^{s,\mu}_{\Psi(t)}-K_N\|_{\mathcal{L}(L^2)}=\sup_{k\in \mathbb{N}:\,k>N} |\sigma_k(t)|.
\end{equation}
Now, using the upper bound in \eqref{eq:cotakappa}:
\begin{equation}\nonumber
\sup_{k>N} |\sigma_k(t)| \le \sup_{k>N}
\frac{C_4}{1+\lambda_k^s[\Psi_a(t)]^\mu} =
\frac{C_4}{1+\lambda_{N+1}^s[\Psi_a(t)]^\mu} \overset{N\to\infty}{\longrightarrow} 0.
\end{equation}
Thus ${_a}\mathcal{K}^{s,\mu}_{\Psi(t)}$ is the uniform limit
 of finite-rank operators. Hence it is compact.
\end{proof}
\begin{remark} Note that for $t=a$, the operator ${_a}\mathcal{K}^{s,\mu}_{\Psi(a)}$ is the identity in $L^2(\Omega)$. 
\end{remark}

 Moreover, the full spectrum is
  \begin{equation}\nonumber
    \sigma\big({}_a\mathcal{K}^{s,\mu}_{\Psi(t)}\big)
    = \sigma_p\big({}_a\mathcal{K}^{s,\mu}_{\Psi(t)}\big) \cup \{0\},
  \end{equation}
  where $0$ is the unique accumulation point; the continuous spectrum is
  $\sigma_c\big({}_a\mathcal{K}^{s,\mu}_{\Psi(t)}\big) = \{0\}$ and the residual spectrum
  is empty, $\sigma_r\big({}_a\mathcal{K}^{s,\mu}_{\Psi(t)}\big) = \varnothing$.

\begin{corollary}[Non-closedness of the range]
    \label{cor:range_not_closed}
    Let $\Omega\subset\mathbb{R}^d$, with $d\le 3$, be a bounded domain and suppose $\frac{d}{4}<s<1$ and $0<\mu_0<\mu_1<1$. For a fixed $t\in(a,T]$ and $\mu\in[\mu_0,\mu_1]$, $\mathrm{Ran}({_a}\mathcal{K}_{\Psi(t)}^{s,\mu})$ is not closed in $L^2(\Omega)$. 
\end{corollary}

\begin{proof}
    Let $\{\varphi_k\}_{k\ge1}$ be the orthonormal basis of eigenfunctions of the Laplacian in $L^2(\Omega)$. For each $n\in\mathbb{N}$, we define:
    \begin{equation} \label{Kvn}
    v^{(n)} \doteq \sum_{k=1}^{n} \varphi_k .
   \end{equation}
   Then
    \begin{equation}\nonumber
    {_a}\mathcal{K}_{\Psi(t)}^{s,\mu}\,v^{(n)}
    =
    \sum_{k=1}^{n} \sigma_k(t)\, \varphi_k .
    \end{equation}
    By Theorem \ref{coro:MittagAsin2}, there exist positive constants $m, M$ depending on $\mu_0, \mu_1$, and $s$ such that, for sufficiently large $k$:
    \begin{equation}\nonumber
    \frac{m}{k^{2s/d}[\Psi_a(t)]^\mu}
    \le
    |\sigma_k(t)|
    \le
    \frac{M}{k^{2s/d}[\Psi_a(t)]^\mu}, \qquad t \in(a,T].
    \end{equation}
    Consequently, $|\sigma_k(t)|^2 \le \frac{M^2}{[\Psi_a(t)]^{2\mu}}k^{-4s/d}$, and since $s>d/4$, we obtain:
    \begin{equation}\nonumber
    \sum_{k=1}^{\infty} |\sigma_k(t)|^2 \leq \frac{M^2}{[\Psi_a(t)]^{2\mu}}\sum_{k=1}^{\infty} k^{-4s/d}<\infty, \qquad t \in (a,T].
    \end{equation}
    Therefore, the function defined as
    \begin{equation}\nonumber
    f \doteq \sum_{k=1}^{\infty} \sigma_k(t)\varphi_k
    \end{equation}
    is in $L^2(\Omega)$ and from (\ref{Kvn}) ${_a}\mathcal{K}_{\Psi(t)}^{s,\mu}(v^{(n)})\overset{n\to\infty}{\longrightarrow}  f$ in $L^2(\Omega)$. This implies that $f\in\overline{\mathrm{Ran}({_a}\mathcal{K}_{\Psi(t)}^{s,\mu})}$. Suppose now that there exists $v\in L^2(\Omega)$ such that ${_a}\mathcal{K}_{\Psi(t)}^{s,\mu}(v)=f$. By considering the spectral representation of $v$
    \begin{equation}\nonumber
    v=\sum_{k=1}^{\infty} v_k\varphi_k,
    \end{equation}
    it follows that:
    \begin{equation}\nonumber
    {_a}\mathcal{K}_{\Psi(t)}^{s,\mu}(v)
    =
    \sum_{k=1}^{\infty}\sigma_k(t)v_k\varphi_k .
    \end{equation}
    Since the right hand side above must be equal to $f$, comparing coefficients it follows that $\sigma_k(t)\, v_k=\sigma_k(t)$ for all $k\in\mathbb{N}$. Since $\sigma_k(t)\neq0$, it follows that $v_k=1$ for all $k\in\mathbb{N}$. This implies:
    \begin{equation}\nonumber
    \|v\|_{L^2(\Omega)}^2
    =
    \sum_{k=1}^{\infty}|v_k|^2
    =
    \sum_{k=1}^{\infty}1
    =
    +\infty,
    \end{equation}
    which contradicts the hypothesis $v\in L^2(\Omega)$. Therefore, $f\notin \mathrm{Ran}({_a}\mathcal{K}_{\Psi(t)}^{s,\mu})$, which proves that range of the operator is not closed in $L^2(\Omega)$.
\end{proof}

\section{Homogeneous Backward Problem with respect to \texorpdfstring{$\Psi$}{Psi}}

% For given $0\le a<T<\infty $, let us consider the following homogeneous forward problem for a bounded set
% $\Omega_T\subset \mathbb{R}^{d+1}$, with a fractional time derivative of order
% $\mu\in(0,1)$ in the Caputo sense with respect to $\Psi$, and of order
% $s\in(0,1)$ in space in the sense of the spectral fractional Laplacian:
% \begin{equation}\label{eq.FBHEDeltaHOMOGPsi}
% (PBH_\Psi)\left\{\begin{aligned} {}_{}^{C}\partial ^{\mu}_{a+;\Psi(t)}
% u(\mathbf{x},t)  &= - (-\Delta_{D,0})^s u(\mathbf{x},t),
% & (\mathbf{x},t) \in \Omega_T, \\
% u(\mathbf{x},t) &= 0, & (\mathbf{x},t) \in \partial\Omega \times (a,T], \\
% u(\mathbf{x},T) &= u_T(\mathbf{x}), & \mathbf{x} \in \Omega.
% \end{aligned}\right.
% \end{equation}

Associated  to problem $(PFH_\Psi)$ (see (\ref{eq.FBHEDeltaPsi}) ) we consider the following two inverse problems:

\begin{problem}
    \label{pro:inverset}
Given the final state $u_T \doteq u_{\mu,s}(\cdot,T)$, find
$u_{\mu,s}(\cdot,t)$ for $a < t < T$.
\end{problem}

\begin{problem}
    \label{pro:inversa}
    Given the final state $u_T \doteq u_{\mu,s}(\cdot,T)$, find the initial state $u_{\mu,s}(\cdot,a^+)$.
\end{problem}

We shall first analyze and solve Problem \ref{pro:inverset}. Let $u_T \doteq
u_{\mu,s}(\cdot,T)\in L^2(\Omega)$ and write it in the form 
\begin{equation}
    \label{eq:usubT}
    u_T(\mathbf{x}) = \sum_{k=1}^\infty (u_T)_k\,\varphi_k(\mathbf{x}),
\end{equation}
where $(u_T)_k=\langle u_T,\varphi_k\rangle_{L^2(\Omega)}$ for
all $k\in \mathbb{N}$. Now using the solution operator (equation \eqref{eq:SolFBHEDeltaPsi}) for the final time $t=T$, we know that 
$u_{\mu,s}(\mathbf{x},T) = {_a}\mathcal{K}_{\Psi(T)}^{s,\mu}(u_a)(\mathbf{x})$,
that is:
\begin{equation}
    \label{eq:coef_uT_Psi}
    u_T(\mathbf{x}) = \sum_{k=1}^\infty \Big[E_{\mu}\!\big(-\lambda_k^{s} (\Psi_a(T))^\mu\big) (u_a)_k\Big] \varphi_k(\mathbf{x}),
\end{equation}
where $(u_a)_k=\langle u_a,\varphi_k\rangle_{L^2(\Omega)}$ for
all $k\in \mathbb{N}$.
Equating the coefficients in \eqref{eq:usubT} and \eqref{eq:coef_uT_Psi}, we
obtain:
\begin{equation}
    (u_T)_k = E_{\mu}\!\big(-\lambda_k^{s} (\Psi_a(T))^\mu\big) (u_a)_k \nonumber
\end{equation}
and therefore
\begin{equation}
    \frac{(u_T)_k}{E_{\mu}\!\big(-\lambda_k^{s} (\Psi_a(T))^\mu\big)} = (u_a)_k, \qquad \forall k \in \mathbb{N}. \label{eq:coefu_T_u_0_Psi}
\end{equation}

Substituting the coefficients \eqref{eq:coefu_T_u_0_Psi} into the general solution given by \eqref{eq:AalphaPsi} and \eqref{eq:SolFBHEDeltaPsi}, for $t\in(a,T]$ we obtain:
\begin{equation}
    u_{\mu,s}(\mathbf{x},t)=
    \sum_{k=1}^\infty \left[ \frac{E_{\mu}\!\big(-\lambda_k^s(\Psi_a(t))^{\mu}\big)}{E_{\mu}\!\big(-\lambda_k^{s} (\Psi_a(T))^\mu\big)}\right](u_T)_k\,   \varphi_k(\mathbf{x}).
\end{equation}
We have then proved the following result: 

\begin{theorem}
    For each $0<a<t\leq T$, there exists a unique solution of Problem \eqref{pro:inverset} that can be writen in the form
\begin{equation}\label{eq.SolFBHEDeltaPsiBack}
u_{\mu,s}(\mathbf{x},t) = \sum_{k=1}^\infty \left[ \
\frac{E_{\mu}\!\big(-\lambda_k^s(\Psi_a(t))^{\mu}\big)}{E_{\mu}\!\big(-\lambda_k^{s}
(\Psi_a(T))^\mu\big)}\right](u_T)_k\varphi_k (\mathbf{x}), \quad
(\mathbf{x},t)\in\Omega_{a,T}.
\end{equation}
\end{theorem}
Note that $u_{\mu,s}(\mathbf{x},t)$  given by (\ref{eq.SolFBHEDeltaPsiBack})  is solution of the following problem:
\begin{equation}\label{eq.FBHEDeltaHOMOGPsi}
(PBH_\Psi)\left\{\begin{aligned} {}_{}^{C}\partial ^{\mu}_{a+;\Psi(t)}
u(\mathbf{x},t)  &= - (-\Delta_{D,0})^s u(\mathbf{x},t),
& (\mathbf{x},t) \in \Omega_T, \\
u(\mathbf{x},t) &= 0, & (\mathbf{x},t) \in \partial\Omega \times (a,T], \\
u(\mathbf{x},T) &= u_T(\mathbf{x}), & \mathbf{x} \in \Omega.
\end{aligned}\right.
\end{equation}

In view of the representation formula \eqref{eq.SolFBHEDeltaPsiBack}, for each
$\mu \in (0,1)$ and $t \in (a,T]$, we define the linear operator
${}_{a}^{T}\mathcal{R}^{\Psi}_{\mu,s}(t) : L^2(\Omega) \longrightarrow
L^2(\Omega)$ by:
\begin{equation}\label{eq:Rmus}
    {}_{a}^{T}\mathcal{R}^{\Psi}_{\mu,s}(t)\, v
    \doteq \sum_{k=1}^{\infty}
    \frac{E_{\mu}\!\big(-\lambda_k^s(\Psi_a(t))^{\mu}\big)}
    {E_{\mu}\!\big(-\lambda_k^s(\Psi_a(T))^{\mu}\big)}
    \, v_k \, \varphi_k,
    \qquad v \in L^2(\Omega),
\end{equation}
where $v_k = \langle v,\varphi_k\rangle_{L^2(\Omega)}$.

In the following theorem we prove several properties of this operator. 

\begin{theorem}  Let $\mu$ be as in Lemma \ref{lem:Mittag}, c).
    For each $t\in (a,T]$, the operator ${}_{a}^{T}\mathcal{R}^{\Psi}_{\mu,s}(t)$ is linear and bounded in $L^2(\Omega)$. Moreover there exists a positive constant $C$ independent of $t$ such that, for all $v\in L^2(\Omega)$ there holds: 
    \begin{equation} \label{cotaR}
    \|{}_{a}^{T}\mathcal{R}^{\Psi}_{\mu,s}(t)v\|_{L^2(\Omega)}
    \le
    C
    \left(
    \frac{\Psi_a(T)}{\Psi_a(t)}
    \right)^{\mu}
    \|v\|_{L^2(\Omega)},\qquad t\in (a,T].
    \end{equation}
\end{theorem}
\begin{proof}
    Let $v\in L^2(\Omega)$ and $t\in(a,T]$. By definition of the operator ${}_{a}^{T}\mathcal{R}^{\Psi}_{\mu,s}(t)$ and the orthonormality of $\{\varphi_k\}_{k\in \mathbb{N}}$ in $L^2(\Omega)$:
    \begin{equation}\label{eq:norma_u_new}
        \|{}_{a}^{T}\mathcal{R}^{\Psi}_{\mu,s}(t)v\|_{L^2(\Omega)}^2
        =
        \sum_{k=1}^{\infty}
        \left|
        \frac{E_\mu\!\big(-\lambda_k^s(\Psi_a(t))^\mu\big)}
        {E_\mu\!\big(-\lambda_k^s(\Psi_a(T))^\mu\big)}
        \right|^2
        |(v,\varphi_k)|_{L^2(\Omega)}^2 .
    \end{equation}

    By Lemma \ref{lem:Mittag}, c) there exist positive constants $C_3, C_4$ such that:
    \begin{equation}\label{abc}
    \frac{E_\mu(-\lambda_k^s(\Psi_a(t))^\mu)}{E_\mu(-\lambda_k^s(\Psi_a(T))^\mu)}
    \le
    \frac{C_4(1+\lambda_k^s(\Psi_a(T))^\mu)}{C_3(1+\lambda_k^s(\Psi_a(t))^\mu)}.
    \end{equation}

    Now, letting $A\doteq\lambda_k^s(\Psi_a(T))^\mu$ and $B\doteq\lambda_k^s(\Psi_a(t))^\mu$ for $k\in \mathbb{N}$ and using that $A\ge B>0$ (since $\Psi_a$ is positive and increasing for $t\in(a,T]$), we obtain the elementary inequality $\frac{1+A}{1+B}\le \frac{A}{B}$. It then follows from (\ref{abc}) that
    \begin{equation}
        \label{eq:estimationRegularization}
            \frac{E_\mu\!\big(-\lambda_k^s(\Psi_a(t))^\mu\big)}
        {E_\mu\!\big(-\lambda_k^s(\Psi_a(T))^\mu\big)}
        \le
        \frac{C_4}{C_3}
        \frac{\lambda_k^s(\Psi_a(T))^\mu}
        {\lambda_k^s(\Psi_a(t))^\mu}
        =:
        C
        \left(
        \frac{\Psi_a(T)}{\Psi_a(t)}
        \right)^{\mu},
    \end{equation}
    where  $C\doteq\dfrac{C_4}{C_3}>0$. Then, from the inequality in \eqref{eq:estimationRegularization} and the identity in \eqref{eq:norma_u_new}, we have:
    \begin{align*}
\|{}_{a}^{T}\mathcal{R}^{\Psi}_{\mu,s}(t)v\|_{L^2(\Omega)}^2
    &\le   \sum_{k=1}^{\infty}
    C^2\left(
    \frac{\Psi_a(T)}{\Psi_a(t)}\right)^{2\mu}
    |(v,\varphi_k)|_{L^2(\Omega)}^2 \\
    &=
    C^2
    \left(
    \frac{\Psi_a(T)}{\Psi_a(t)}
    \right)^{2\mu}
    \|v\|_{L^2(\Omega)}^2 ,\qquad t\in(a,T]
    \end{align*}
    which proves the boundedness of the operator ${}_{a}^{T}\mathcal{R}^{\Psi}_{\mu,s}(t)$ for all $t\in(a,T]$ and hence the stability of the backward inverse problem \ref{pro:inverset}. This completes the proof.
\end{proof}

By solving the backward inverse Problem \ref{pro:inverset}, we were able to write the solution
at times $t>a$ in the form
\begin{equation}
  u_{\mu,s}(\mathbf{x},t)={}_{a}^{T}\mathcal{R}^{\Psi}_{\mu,s}(t)\,u_T(\mathbf{x}),
\qquad (\mathbf{x},t) \in \Omega_{a,T},
\end{equation}
where ${}_{a}^{T}\mathcal{R}^{\Psi}_{\mu,s}(t)$ is the operator defined in (\ref{eq:Rmus}). Although, as shown, this operators are bounded for all $t>a$, 
%
%
%
% \begin{remark}[Continuous dependence on final data]
%     Let $u_T, \tilde u_T \in L^2(\Omega)$ and let $t\in(a,T]$ be fixed. Denote by $u_{\mu,s}(\cdot,t)$ and $\tilde u_{\mu,s}(\cdot,t)$ the solutions of the inverse problem $P_t^\Psi$ associated with the final data $u_T$ and $\tilde u_T$, respectively. By the linearity of the operator ${}_{a}^{T}\mathcal{R}^{\Psi}_{\mu,s}(t)$:
%     \begin{equation}\nonumber
%     u_{\mu,s}(\cdot,t)-\tilde u_{\mu,s}(\cdot,t)
%     =
%     {}_{a}^{T}\mathcal{R}^{\Psi}_{\mu,s}(t)(u_T-\tilde u_T).
%     \end{equation}
%     Using the stability estimate of the operator:
%     \begin{equation}\nonumber
%     \|u_{\mu,s}(\cdot,t)-\tilde u_{\mu,s}(\cdot,t)\|_{L^2(\Omega)}
%     \le
%     C
%     \left(
%     \frac{\Psi_a(T)}{\Psi_a(t)}
%     \right)^{\mu}
%     \|u_T-\tilde u_T\|_{L^2(\Omega)}.
%     \end{equation}
%     In particular, for each $t\in(a,T]$, the mapping $u_T \longmapsto u_{\mu,s}(\cdot,t)$ is continuous from $L^2(\Omega)$ into itself.
% \end{remark}
the factor $\left( \frac{\Psi_a(T)}{\Psi_a(t)}
\right)^{\mu}$ in the right hand side of (\ref{cotaR}) becomes unbounded as $t\to a^+$, reflecting the ill-posed nature
of the fractional backward inverse Problem \ref{pro:inversa}. Indeed, if
$t\to a^+$, then $\Psi_a(t)\to 0$ and consequently
$E_{\mu}(-\lambda_k^s(\Psi_a(t))\to 1$. Therefore, letting $t\to a^+$ in (\ref{eq.SolFBHEDeltaPsiBack}) we formally obtain:
\begin{equation}
    \label{eq:SolBackInitial}
u_{\mu,s}(\cdot,a^+) = \sum_{k=1}^\infty \frac{\langle u_T,\varphi_k\rangle}
{E_{\mu}\!\big(-\lambda_k^s(\Psi_a(T))^{\mu}\big)} \,\varphi_k.
\end{equation}
It is clear that his expresion is only formal since conditions on $u_T$ a required for the RHS of (\ref{eq:SolBackInitial}) to define a function in $L^2(\Omega)$. These are precisely the so called Picard's conditions (see \cite{Engl1996}). Namely, that 
\begin{equation}
\sum_{k=1}^\infty \left|\frac{\langle u_T,\varphi_k\rangle}
{E_{\mu}\!\big(-\lambda_k^s(\Psi_a(T))^{\mu}\big)}\right|^2\;<\;\infty  
\end{equation}
Moreover, the expresion in the RHS of (\ref{eq:SolBackInitial}) is precisely $\left({_a}\mathcal{K}_{\Psi(T)}^{s,\mu}\right)^\dagger$, the Moore-Penrose generalized inverse of the operator $ {_a}\mathcal{K}_{\Psi(T)}^{s,\mu}$ defined in (\ref{eq:AalphaPsi}), for $t=T$. 

\medskip
We shall show now that this operator is unbounded, what leads precisely to the ill-posedness of Problem \ref{pro:inversa}. For that assume $u_T\in Ran({_a}\mathcal{K}_{\Psi(T)}^{s,\mu})$ and let $\tilde u_T\in L^2(\Omega)$ be a perturbation of this final data with
$\tilde u_{\mu,s}$ its corresponding solution. Then:
\begin{equation}\nonumber
u_{\mu,s}(\cdot,a^+)-\tilde u_{\mu,s}(\cdot,a^+) = \sum_{k=1}^\infty
\frac{(u_T-\tilde u_T,\varphi_k)}
{E_{\mu}\!\big(-\lambda_k^s(\Psi_a(T))^{\mu}\big)} \,\varphi_k.
\end{equation}
Applying Parseval's identity it follows that
\begin{equation}
    \label{eq.ParsevalPsiIll}
    \|u_{\mu,s}(\cdot,a^+)-\tilde u_{\mu,s}(\cdot,a^+)\|_{L^2(\Omega)}^2
    =
    \sum_{k=1}^\infty
    \frac{|\langle u_T-\tilde u_T,\varphi_k\rangle|^2}
    {E_{\mu}\!\big(-\lambda_k^s(\Psi_a(T))^{\mu}\big)^2}.
\end{equation}
Now, from the upper bound in Lemma \ref{lem:Mittag} c) for $t=T$, we have:
\begin{equation}\nonumber
\frac{1}{E_{\mu}\!\big(-\lambda_k^s(\Psi_a(T))^{\mu}\big)} \ge
\frac{1+\lambda_k^s(\Psi_a(T))^{\mu}}{C_4}>0, \quad \forall k\in\mathbb{N}.
\end{equation}
Substituting into \eqref{eq.ParsevalPsiIll} we obtain
\begin{equation}
    \label{eq.IllPsi}
    \|u_{\mu,s}(\cdot,a^+)-\tilde u_{\mu,s}(\cdot,a^+)\|_{L^2(\Omega)}^2
    \ge
    \sum_{k=1}^\infty
    \frac{\big(1+\lambda_k^s(\Psi_a(T))^{\mu}\big)^2}{C_4^2}
    |\langle u_T-\tilde u_T,\varphi_k\rangle|^2.
\end{equation}
Now, given $\delta>0$, define $\tilde u_{T,k} = u_T + \delta\,\varphi_k$
for $k\in\mathbb{N}$, and let $\tilde u_{\mu,s,k}$ be its corresponding solution. Then $\|u_T-\tilde u_{T,k}\|_{L^2(\Omega)}=\delta$ for all $k\in\mathbb{N}$ and
from \eqref{eq.IllPsi}, with  it follows that
\begin{equation}\nonumber
\|u_{\mu,s}(\cdot,a^+)-\tilde u_{\mu,s,k}(\cdot,a^+)\|_{L^2(\Omega)} \ge \delta\; 
\frac{1+\lambda_k^s(\Psi_a(T))^{\mu}}{C_4}\to +\infty, \quad \text{as }
k\to\infty.
\end{equation}
This proves the operator $\left({_a}\mathcal{K}_{\Psi(T)}^{s,\mu}\right)^\dagger$ is unbounded.

In general, the previous result can be interpreted within the framework of
linear inverse problems in Hilbert spaces. Indeed, problem $(PBH_{\Psi})$ can
be reformulated as a linear problem in $L^2(\Omega)$. The relationship
between the initial state $u_a$ and the final state $u_T$ is given by:
\begin{equation}
    \label{eq:BackwardOperator}
    {_a}\mathcal{K}_{\Psi(T)}^{s,\mu}(u_a)=u_T.
\end{equation}
As seen, the operator ${_a}\mathcal{K}_{\Psi(T)}^{s,\mu}$ is self-adjoint,
non-negative, injective, and compact (see Proposition~\ref{prop:Operator_Sol}
for $t=T$). Hence, its singular system coincides with its eigensystem and it is given by
$(\sigma_k(T),\varphi_k,\varphi_k)_{k\in\mathbb{N}}$, where $\sigma_k(T) =
E_{\mu}\!\big(-\lambda_k^s(\Psi_a(T))^\mu\big)$.

\medskip
According to the classical theory of linear inverse problems
(cf.~\cite{Engl1996}), the Moore-Penrose generalized inverse
$({_a}\mathcal{K}_{\Psi(T)}^{s,\mu})^\dagger$, which in this case coincides with the clasical inverse due to injectivity, is a densely defined, unbounded
linear operator with a closed graph. Hence, the inverse problem associated to \eqref{eq:BackwardOperator} is ill-posed. 

% Assuming $u_T \in \mathrm{Ran}(\mathcal{K}^{\Psi(T)}_{s,\mu})$, then by
% Picard's criterion, $(\mathcal{K}^{\Psi(T)}_{s,\mu})^\dagger u_T \in
% L^2(\Omega)$ and:
% \begin{equation}\nonumber
% u^{\dagger}_{\mu,s}(\cdot,a) = (\mathcal{K}^{\Psi(T)}_{s,\mu})^\dagger
% u_T=\sum_{k=1}^\infty \frac{\langle u_T,\varphi_k\rangle}
% {E_{\mu}\!\big(-\lambda_k^s(\Psi_a(T))^{\mu}\big)} \,\varphi_k.
% \end{equation}
% Comparing with \eqref{eq:SolBackInitial}, this exactly coincides with the
% spectral expansion obtained for the backward solution of $(PFH_\Psi)$, i.e.,
% $u_{\mu,s}(\cdot,a^+) = (\mathcal{K}^{\Psi(T)}_{s,\mu})^\dagger u_T$.
% Therefore, $u_{\mu,s}(\cdot,a^+)$ constitutes the unique solution
% in $L^2(\Omega)$ for the inverse problem \eqref{eq:BackwardOperator}.

%%%%%%%%%%%%%%     H A S T A 
\medspace
\begin{remark}[Interpretation as a regularization family]
    For each $\mu \in (0,1)$, the family of operators $\{ {}_{a}^{T}\mathcal{R}^{\Psi}_{\mu,s}(t) \}_{t\in(a,T)}$, as defined in (\ref{eq:Rmus}) allows for the construction of a \textit{regularization family} associated for the backward problem $(PBH_{\Psi})$. We define the family of operators $\{ \mathcal{R}_\alpha \}_{\alpha \in (0, T-a)}$ by:
    \begin{equation}
        \mathcal{R}_\alpha \doteq {}_{a}^{T}\mathcal{R}^{\Psi}_{\mu,s}(a+\alpha), \quad \text{where } \alpha \doteq t-a > 0.
    \end{equation}
    For each $\alpha$, the operator $\mathcal{R}_\alpha$ acts on the final data $u_T$ via spectral expansion. It turns out that these operators constitute a regularization family for $({_a}\mathcal{K}_{\Psi(T)}^{s,\mu})^\dagger$. In fact: 
    \begin{enumerate}
        \item Stability: For each fixed $\alpha \in (0, T-a)$, $\mathcal{R}_\alpha$ is a bounded linear operator. The loss of boundedness occurs only at the critical limit $\alpha \to 0^+$.
        \item Convergence: 
        $\mathcal{R}_\alpha u_T \overset{\alpha\to 0}{\longrightarrow} ({_a}\mathcal{K}_{\Psi(T)}^{s,\mu})^\dagger\, u_T$, pointwise for  $u_T \in \mathcal{D}(({_a}\mathcal{K}_{\Psi(T)}^{s,\mu})^\dagger).$
    \end{enumerate}
    Note that the backward diffusion operator at $t=a+\alpha > a$ acts as a frequency filter that stabilizes the inversion, mitigating the noise amplification at high-frequency modes.
\end{remark}

\subsection{Stability Analysis and Relationship between Problems}

Problems \ref{pro:inverset}  and \ref{pro:inversa} have very different properties, mainly regarding their stability.  In fact:
\begin{enumerate}
    \item \textbf{Problem \ref{pro:inversa} is ill-posed:} recovering the exact initial state $u_{\mu,s}(\cdot, a^+)$ is unstable. This is so because the forward operator is compact and it has non-closed range. Hence its Moore-Penrose inverse is unbounded (not continuous). Noise at high frequencies in the data amplifies unboundedly when reversing the process to $t=a$.
    \item \textbf{Problem \ref{pro:inverset} is well-posed:} the problem of finding an intermediate state $u_{\mu,s}(\cdot, t)$ for $t \in (a, T)$ is \textit{well-posed}. The inverse solution operator remains bounded because the Mittag-Leffler kernel acts as a filter which is singular only at $t=a$.
\end{enumerate}

\subsection{Degree of Ill-posedness}
The degree of ill-posedness of an operator equation $Kx=y$ is characterized by the decay rate of its singular values $\{\sigma_k\}$ (see \cite{Engl1996}). The problem is said to be mildly or moderately ill-posed if $\sigma_k = O(k^{-\alpha})$ for some $\alpha > 0$, and severely ill-posed otherwise, e.g. if such a decay rate is exponential, that is $\sigma_k = O(e^{-k})$.
Thus, in the classical case ($a=0$, $\Psi(t)=t$, $\mu=s=1$), the problem is \textit{severely
ill-posed} since, in that case,  the eigenvalues of the operator ${_0}\mathcal{K}_{t}^{1,1}$ are given by  $e^{-\lambda_k t }$ (see (\ref{eq.SolFBHEDeltaCaputo3})). In
contrast, under fractional diffusion ($0 < \mu < 1$), the behavior is
significantly more stable:
\begin{enumerate}
    \item \textbf{Stability in Problem \ref{pro:inverset} ($t > a$):} The high-frequency amplification factors are bounded. In fact as it can be seen in (\ref{eq:Rmus}), for any $a<t\le T$, the eigenvalues of ${}_{a}^{T}\mathcal{R}^{\Psi}_{\mu,s}(t)$ satisfy $\frac{\sigma_k(t)}{\sigma_k(T)}=\frac{E_{\mu}\big(-\lambda_k^s(\Psi_a(t))^\mu\big)}{E_{\mu}\big(-\lambda_k^s(\Psi_a(T))^\mu\big)} \to \left( \frac{\Psi_a(T)}{\Psi_a(t)} \right)^\mu\;<\infty$, as $k\to\infty$. Thus, the  operator ${}_{a}^{T}\mathcal{R}^{\Psi}_{\mu,s}(t)$, asociated to the inversion of problem ($PBH_\Psi$) in (\ref{eq.FBHEDeltaHOMOGPsi}) is bounded, making Problem \ref{pro:inverset} \textbf{well-posed}.
    \item \textbf{Regularity in Problem \ref{pro:inversa} ($t \to a^+$):} From Theorem \ref{coro:MittagAsin2}, it follows that the eigenvalues $\sigma_k$ of the operator ${_a}\mathcal{K}_{\Psi(t)}^{s,\mu}$ defined in (\ref{eq:AalphaPsi}) clearly satisfy $\sigma_k=O(k^{-2s/d})$  and therefore Problem \ref{pro:inversa} is only    \textbf{moderately ill-posed}.
\end{enumerate}

\begin{theorem}[Uniform limit of the modal factor]
    Let $0\le a<a_0\le T$, $0<s_0\le1$, and $0<\mu_0<\mu_1<1$. For all $t\in[a_0,T]$, $\lim_{k\to\infty}\frac{\sigma_k(t)}{\sigma_k(T)} = \big(\frac{\Psi_a(T)}{\Psi_a(t)}\big)^\mu$, and the convergence is uniform with respect to $(t,s,\mu)\in[a_0,T]\times[s_0,1]\times[\mu_0,\mu_1]$.
\end{theorem}

\begin{proof}
    From the asymptotic expansion of the Mittag-Leffler function for $p=1$ and using \eqref{eq:argumentoZk}, we obtain:
    \begin{equation}\nonumber
    \sigma_k(t)
    =
    \frac{1}{\Gamma(1-\mu)}
    \lambda_k^{-s}
    [\Psi_a(t)]^{-\mu}
    +
    \mathcal O\!\left(
    \lambda_k^{-2s}[\Psi_a(t)]^{-2\mu}
    \right),
    \qquad k\to\infty .
    \end{equation}

    Factoring out $\lambda_k^{-s}$, we can write:
    \begin{equation}\nonumber
    \sigma_k(t)
    =
    \lambda_k^{-s}
    \left(
    \frac{[\Psi_a(t)]^{-\mu}}{\Gamma(1-\mu)}
    +
    r_k(t)
    \right),
    \end{equation}
    where there exists a constant $C>0$ such that:
    \begin{equation}\nonumber
    |r_k(t)|
    \le
    C\,\lambda_k^{-s}[\Psi_a(t)]^{-2\mu}.
    \end{equation}

    If $t\in[a_0,T]$ and $s\ge s_0$, since $\Psi$ is strictly increasing, we have $\Psi_a(t)\ge\Psi_a(a_0)>0$, and therefore:
    \begin{equation}\nonumber
    |r_k(t)|
    \le
    C\,\lambda_k^{-s_0}[\Psi_a(a_0)]^{-2\mu}.
    \end{equation}

    Furthermore, since $\mu\in[\mu_0,\mu_1]$, it follows that:
    \begin{equation}\nonumber
    [\Psi_a(a_0)]^{-2\mu}
    \le
    [\Psi_a(a_0)]^{-2\mu_1}.
    \end{equation}

    Thus,
    \begin{equation}\nonumber
    \sup_{t\in[a_0,T]}
    |r_k(t)|
    \le
    C\,\lambda_k^{-s_0}[\Psi_a(a_0)]^{-2\mu_1}
    \longrightarrow 0
    \qquad (k\to\infty),
    \end{equation}
    which implies that $r_k(t)\to 0$ uniformly in $(t,s,\mu)\in[a_0,T]\times[s_0,1]\times[\mu_0,\mu_1]$.

    Denoting $c_{\mu}=\frac{1}{\Gamma(1-\mu)}$, we can write:
    \begin{equation}\nonumber
    \frac{\sigma_k(t)}{\sigma_k(T)}
    =
    \frac{
        c_{\mu}[\Psi_a(t)]^{-\mu}+r_k(t)
    }{
        c_{\mu}[\Psi_a(T)]^{-\mu}+r_k(T)
    },
    \end{equation}
    and consequently, we obtain:
    \begin{equation}\nonumber
    \frac{\sigma_k(t)}{\sigma_k(T)}- \left(\frac{\Psi_a(T)}{\Psi_a(t)}\right)^{\mu} = \frac{ [\Psi_a(T)]^{-\mu} r_k(t)-[\Psi_a(t)]^{-\mu}r_k(T) }{ [\Psi_a(T)]^{-\mu}(c_{\mu}[\Psi_a(T)]^{-\mu}+r_k(T)) }.
    \end{equation}
    Since $r_k(T)\to 0$ uniformly with respect to $s$ and $\mu$, there exists $k_0$ such that for $k\ge k_0$:
    \begin{equation}\nonumber
    \left|
    c_{\mu}[\Psi_a(T)]^{-\mu}+r_k(T)
    \right|
    \ge
    \frac12\,c_{\mu}[\Psi_a(T)]^{-\mu}.
    \end{equation}

    As a result,
    \begin{equation}\nonumber
    \left|
    \frac{\sigma_k(t)}{\sigma_k(T)}
    -
    \left(\frac{\Psi_a(T)}{\Psi_a(t)}\right)^\mu
    \right|
    \le
    C\big(|r_k(t)|+|r_k(T)|\big).
    \end{equation}

    Given $\varepsilon>0$, since $r_k(t)\to 0$ uniformly in $(t,s,\mu)$, there exists $K$ such that for $k\ge K$:
    \begin{equation}\nonumber
    |r_k(t)|<\frac{\varepsilon}{4C},
    \qquad
    |r_k(T)|<\frac{\varepsilon}{4C}.
    \end{equation}

    Therefore,
    \begin{equation}\nonumber
    \left|
    \frac{\sigma_k(t)}{\sigma_k(T)}
    -
    \left(\frac{\Psi_a(T)}{\Psi_a(t)}\right)^\mu
    \right|
    <
    \varepsilon,
    \end{equation}
    uniformly with respect to $(t,s,\mu)\in[a_0,T]\times[s_0,1]\times[\mu_0,\mu_1]$, which concludes the proof.
\end{proof}

In this sense, the $\Psi$-fractional backward problem constitutes an example of
an inverse problem whose severity is strictly lower than that of \emph{mildly
ill-posed} problems in the traditional classification.
\section{Conclusions}
In this work, we studied a generalized space-time fractional heat equation defined on an open and bounded domain $\Omega \subset \mathbb{R}^d$. The model incorporates a time-fractional derivative of order $\mu \in (0,1)$ in the Caputo sense with respect to a time scale function $\Psi$, alongside a Dirichlet spectral fractional Laplacian $(-\Delta_{D,0})^s$ of order $s \in (0,1)$ in the spatial variable.

The main contributions and findings of this study are summarized as follows:
\renewcommand{\theenumi}{\roman{enumi}}
\begin{enumerate}
    \item \textbf{Analysis of the Direct Problem:} By applying a spectral decomposition based on the orthonormal eigenfunctions $\{\varphi_k\}_{k\in\mathbb{N}}$ of the Dirichlet Laplacian, we established the existence, uniqueness, and explicit representation of the solution to the forward problem $(PFH_\Psi)$. The solution operator ${}_a\mathcal{K}_{\Psi(t)}^{s,\mu}$, expressed via the two-parameter Mittag-Leffler function $E_\mu$, encompasses several well-known fractional diffusion models as special cases, including the classical Caputo ($\Psi(t)=t$), Caputo-Hadamard ($\Psi(t)=\ln t$), and Caputo-Katugampola ($\Psi(t)=t^\rho/\rho$) frameworks.

    \item \textbf{Properties of the Solution Operator:} We rigorously analyzed the spectral behavior and functional properties of ${}_a\mathcal{K}_{\Psi(t)}^{s,\mu}$. It was proven that for any intermediate time $t \in (a,T]$, the solution operator is self-adjoint, positive, compact, and injective in $L^2(\Omega)$. Furthermore, for $s > d/4$, we proved that the range $\mathrm{Ran}({}_a\mathcal{K}_{\Psi(t)}^{s,\mu})$ is not closed, directly confirming the inherent ill-posedness of the corresponding inverse problems.

    \item \textbf{Backward Inverse Problems and Stability Hierarchy:} We formulated and investigated two distinct backward inverse problems from final observational data $u_T = u(\cdot, T)$:
    {\renewcommand{\theenumi}{\alph{enumi}}
    \begin{enumerate}
        \item \textbf{Intermediate State Reconstruction (Problem 4.1):} We demonstrated that finding the state $u(\cdot, t)$ at any intermediate time $t \in (a, T)$ is a \textbf{well-posed problem}. The inverse solution operator ${}_a^T\mathcal{R}_{\mu,s}^\Psi(t)$ is linear and bounded in $L^2(\Omega)$, as the decaying Mittag-Leffler kernel effectively acts as a frequency filter that suppresses high-mode instability.
        \item \textbf{Initial State Reconstruction (Problem 4.2):} Reconstructing the initial state $u(\cdot, a^+)$ is inherently \textbf{ill-posed} due to the unboundedness of the generalized Moore-Penrose inverse $({}_a\mathcal{K}_{\Psi(T)}^{s,\mu})^\dagger$. However, the family of intermediate backward operators $\{{}_a^T\mathcal{R}_{\mu,s}^\Psi(a+\alpha)\}_{\alpha > 0}$ forms a natural, temporal regularization scheme for the inverse problem as $\alpha \to 0^+$.
    \end{enumerate}
    }

    \item \textbf{Degree of Ill-Posedness and Asymptotic Behavior:} Utilizing the spectral $s$-Weyl estimate and the asymptotic behavior of the Mittag-Leffler function, we established that the singular values $\sigma_k(T)$ decay polynomially as $\mathcal{O}(k^{-2s/d})$. This demonstrates that the $\Psi$-Caputo fractional backward problem is only \textbf{moderately ill-posed}, in sharp contrast to the classical backward heat equation ($a=0$, $\Psi(t)=t$, $\mu = s = 1$), which exhibits an exponential singular value decay $\mathcal{O}(e^{-k})$ and is severely ill-posed. Additionally, we proved that the mode amplification factor $\sigma_k(t)/\sigma_k(T)$ converges uniformly to $(\Psi_a(T)/\Psi_a(t))^\mu$ as $k \to \infty$, showing how the choice of the temporal scale function $\Psi$ directly influences and controls spectral mode amplifications.
\end{enumerate}

In conclusion, replacing the physical time variable $t$ with a scaling function $\Psi(t)$ provides a versatile and mathematically robust framework for modeling subdiffusive and ultra-slow processes in bounded domains, while significantly softening the intrinsic instability associated with time-reversed historical reconstruction. Future research may focus on numerical implementations, parameter identification for unknown orders $\mu$ or $s$, and the analysis of non-homogeneous source terms under noisy data.

\bibliographystyle{plain}
\bibliography{reference2}
\end{document}